\documentclass[a4paper,10pt]{article}

\usepackage[latin1]{inputenc}
\usepackage{natbib}

\usepackage{amsmath}
\usepackage{amssymb}
\usepackage{amsthm}
\usepackage{bbm}
\usepackage{enumerate}
\usepackage{graphics}

\newtheorem{theorem}{Theorem}[section]
\newtheorem{lemma}[theorem]{Lemma}

\newtheorem{hypothesis}{Hypothesis}
\newtheorem{proposition}[theorem]{Proposition}

\newtheorem{remark}[theorem]{Remark}
\newtheorem{corollary}[theorem]{Corollary}

\newcommand{\E}{\mathbb{E}}

\newcommand{\p}{\mathbb{P}}
\newcommand{\R}{\mathbb{R}}
\newcommand{\psimax}{\psi_{\max}}

\newcommand{\osc}{\operatorname{osc}}
\newcommand{\ue}{U^{\epsilon}}
\newcommand{\uec}{U^{\epsilon,c}}
\newcommand{\uei}{U^{\epsilon,\delta i, \delta (i+1)}}

\newcommand{\bt}{\bar{\theta}}
\newcommand{\Ep}{\E_{N}}
\newcommand{\pp}{\p_{N}}

\newcommand{\TV}{\mathrm{TV}}
\renewcommand{\epsilon}{\varepsilon}
\renewcommand{\mid}{\,|\,}

\def \ds{\displaystyle}

\def\build#1#2#3{\mathrel{\mathop{\kern 0pt#1}\limits_{#2}^{#3}}}

\renewcommand{\thefootnote}{\textit{\alph{footnote}}}

\begin{document}

\hspace*{2cm}
\begin{center}
{\LARGE Fast simulated annealing in $\R^d$ and an}\\[2mm]
{\LARGE application to maximum likelihood estimation}\\[2mm]
{\LARGE in state-space models}

\vspace{4mm}
{\large Sylvain RUBENTHALER\footnote{Universit\'e de Nice - Sophia Antipolis,
  Laboratoire Dieudonn\'e, Parc Valrose, 06108 Nice C\'edex 02, France,
  {\tt rubentha@math.unice.fr }}},\\[2mm]
{\large Tobias RYD\'EN\footnote{Centre for Mathematical Sciences,
  Lund University, Box~118, 221~00 Lund, Sweden,
  {\tt tobias.ryden@matstat.lu.se} or {\tt magnus.wiktorsson@matstat.lu.se}.
  Both authors were supported by grants from the Swedish National
  Research Council.}
and Magnus WIKTORSSON$^{\thefootnote}$}

\vspace{3mm}
{\large 12 September 2006}
\end{center}

\vspace{3mm}

\begin{abstract}
Using classical simulated annealing to maximise a function
$\psi$ defined on a subset of $\R^d$, the probability
$\p(\psi(\theta_n)\leq \psi_{\max}-\epsilon)$ tends to zero
at a logarithmic rate as $n$ increases;
here $\theta_n$ is the state in the $n$-th stage of the simulated
annealing algorithm and $\psi_{\max}$ is the maximal value of
$\psi$. We propose a modified scheme for which this probability
is of order $n^{-1/3}\log n$, and hence vanishes at an algebraic rate.
To obtain this faster rate, the exponentially decaying acceptance
probability of classical simulated annealing is replaced by a more
heavy-tailed function, and the system is cooled faster.
We also show how the algorithm may be applied to functions that
cannot be computed exactly but only approximated, and give
an example of maximising the log-likelihood function for a
state-space model.
\end{abstract}

\noindent
{\it Keywords~:}  Central limit and other weak theorems, Computational methods in Markov chains, Sequential estimation, Markov processes with continuous parameter, Monte Carlo methods, Stochastic programming.

\noindent
{\it MSC~:} 60F05, 60J22, 60J25, 62L12,  65C05, 82C80, 90C15.

\section{Introduction}\label{sec:Intro}

Simulated annealing is a simulation-based approach to the problem
of optimising a function.
In the present paper we will be concerned
with a real-valued function, $\psi$ say, defined on a subset
$\Theta$ of $\R^d$, and our aim is to \textit{maximise} $\psi$.
Thus we assume that $\psi$ is bounded and that its supremum is
attained at least at one point. Simulated
annealing is designed to find the global maximum of $\psi$, even if
$\psi$ has local maxima. It has been
extensively studied, see for instance
\cite{delmoral-miclo-99}, \cite{catoni-99} and \cite{catoni-cot-98}
among many others,
and \cite{bartoli-delmoral-01} for an elementary introduction to the
subject. The classical simulated annealing algorithm departs from
a Markov transition kernel, which we denote by $K(\cdot,\cdot)$,
on $\Theta$, and a positive sequence $(\beta_n)_{n\geq 0}$ increasing to
infinity. The sequence $(\beta_n)$ is often referred to as an
(inverse) \textit{cooling schedule}, because $1/\beta_n$ is often
interpreted as a temperature; this terminology originates from
statistical physics.
Then, starting from an initial point $\theta_0\in\Theta$,
a sequence $(\theta_n)_{n\geq 0}$ is constructed recursively as
follows.
\begin{itemize}
\item[(a1)] In stage $n$, given the current state $\theta_n$, sample
  a new proposed position $Z$ from $K(\theta_n,\cdot)$.
\item[(a2)] Set $\theta_{n+1}=Z$ with probability
  $$\exp(-\beta_n (\psi(\theta_n)-\psi(Z))_+)$$
  and $\theta_{n+1}=\theta_n$ otherwise.
\end{itemize}
Here $(\cdot)_+$ is the positive part. We notice that if
$\psi(Z)\geq \psi(\theta_n)$, then the proposed new state $Z$ is accepted
with probability one. A proposal $Z$ at which $\psi$ is smaller than
at the current $\theta_n$ may be accepted, but this becomes increasingly
unlikely for large $n$ since $\beta_n\to\infty$.

The basic idea of simulated annealing is as follows. The update rule
above corresponds to a Markov transition kernel, $K_\beta$ say,
on $\Theta$; cf.\ (\ref{Eq:defK}) below.
Under additional assumptions including that
$K$ is positive recurrent and reversible with respect to its stationary
distribution, $\gamma$ say,
$$
\gamma(dx)\,K(x,dy)=\gamma(dy)\,K(y,dx),
$$
one can prove that for fixed $\beta$, the stationary distribution of $K_\beta$
is absolutely continuous with respect to $\gamma$
with Radon-Nikodym derivative proportional to $\exp\{\beta \psi(x)\}$
(cf. \citealp{catoni-99}, Proposition~1.2, or
\citealp{bartoli-delmoral-01}, p.~64).
This indicates that as $\beta$ increases, this stationary distribution
becomes increasingly concentrated around the maxima of $\psi$.
Now, in the beginning of the simulation scheme
$\beta_n$ is small (the temperature is high), and the particle $\theta_n$
is allowed to explore the space $\Theta$ rather freely. When the temperature
cools down ($\beta_n$ gets large), the particle is more and more lured
to the regions where
$\psi$ is large and should in the limit end up at a maximum point of $\psi$.

Obviously, the kernel $K$ and the sequence $(\beta_n)$ are important
design parameters of the algorithm. A typical choice for $(\beta_n)$
is a logarithmic increase; $\beta_n=\beta_0\log(n+e)$ for some $\beta_0>0$.
We note that with this cooling schedule, the acceptance probability
in (a2) above becomes
\begin{equation}\label{eq:stdaccprob}
  (n+e)^{-\beta_0 (\psi(\theta_n)-\psi(Z))_+ }.
\end{equation}

Under additional regularity assumptions one can prove that for
$\beta_0$ small enough and if $\psi$ has a single global maximum,
it holds that for all $\epsilon>0$,
\begin{equation}\label{eq:simannconsistency}
\p(\psi(\theta_n)\leq \psimax -\epsilon) \to 0 \quad\mbox{as $n\to\infty$},
\end{equation}
where $\psimax=\sup_{x\in\Theta}\psi(x)$. How fast is this convergence?
In many works on simulated annealing the space $\Theta$ is assumed
finite, and one may then let $\epsilon\to 0$ and thus study
$\p(\psi(\theta_n)< \psimax)$. Typically this probability tends to
zero at an algebraic rate, see for instance
\citet[Eq.~(22)]{gielis-maes-99} (take $f$ as the indicator function
of non-optimal states) and references in this
paper. For a continuous $\Theta$ the situation is different.
If $\Theta\subset\R^d$ one can show
(see Appendix~\ref{app:classicalrate}) that the rate of convergence in
(\ref{eq:simannconsistency}) is only logarithmic. Alternatively,
one can prove that there are numbers $C$ and $C'$ such that
for any $\epsilon>0$,
\begin{equation}\label{eq:stdrateconv}
\p(\psi(\theta_n)\leq \psimax -\epsilon) \leq C n^{-C'\epsilon}
\end{equation}
Thus, the algebraic rate becomes infinitely slow as $\epsilon\to 0$.
\citet{locatelli-01} proposed a refinement of the annealing scheme
that reaches non-vanishing algebraic rates,
but it requires knowledge of $\psimax$ which is an assumption
we do not want to make.

In the present paper we propose a modified simulated annealing scheme
such that for any $\epsilon>0$ there is a number $C_\epsilon$ such that
\begin{equation}\label{eq:rate13conv}
\p\left(\psi(\theta_n)\leq \psi_{\text{max}} -\epsilon\right)\leq
  C_\epsilon n^{-1/3}(1+\log n).
\end{equation}
We will then say that the rate of convergence is $1/3$, up to a
logarithmic term.

\section{Description of the new simulated annealing scheme}
\label{sec:description}

Just as in classical simulated annealing, the proposed scheme departs
from a Markov transition kernel $K$ and a cooling schedule $(\beta_n)$.
The difference lies in that the exponential function of the
classical algorithm's update is replaced by a different function, and that
the cooling schedule is altered. More precisely, we let
$g:\R^+\to \R^+$ be a $C^\infty$-function such that $g(0)=1$,
$g$ is non-decreasing and $g(t)\to \infty$ as $t\to\infty$.
We set $f=1/g$ and suppose that $f$ is convex and such that
$\sup_{t\geq 0}|t f'(t)|<\infty$. Then the algorithm looks as follows.
\begin{itemize}
\item[(b1)] In stage $n$, given the current state $\theta_n$,
  sample a new proposed position $Z$ from $K(\theta_n,\cdot)$.
\item[(b2)] Set $\theta_{n+1}=Z$ with probability
  $$
      f(\beta_n(\psi(\theta_n)-\psi(Z))_+)
  $$
  and $\theta_{n+1}=\theta_n$ otherwise.
\end{itemize}
In classical simulated annealing $g(t)=e^t$.
In Section~\ref{sec:rate} we advocate the particular choice
$g(t)=1+t/\tau$ for some $\tau>0$, and thus $f(t)\sim \tau/t$ as $t\to\infty$.
Compared to $f(t)\sim\exp(-t)$, this allows the algorithm to be `more bold'
in exploring regions far away from the current state. On the other hand
we will let $\beta_n$ be of order $n^\alpha$ with $\alpha=1/3$,
so that this sequence increases much faster than logarithmically.
Together, these conditions imply (\ref{eq:rate13conv}).
We also remark that with $g$ as above and $\beta_n = n^{1/3}$,
the acceptance probability in (b2) becomes
\begin{equation}\label{eq:newaccprob}
  \frac{1}{1+ \frac{n^{1/3}}{\tau}(\psi(\theta_n)-\psi(Z))_+},
\end{equation}
which should be compared to (\ref{eq:stdaccprob}); we see that
(\ref{eq:newaccprob}) decays much slower as
$\psi(\theta_n)-\psi(Z)\to\infty$, and thus again that the new
algorithm is less likely to reject proposals with function values
far below the current one.

Modifications of the acceptance function $f(t)=\exp(-t)$ of
classical simulated annealing to speed
up convergence rates have been discussed extensively in the
statistical physics literature, and is there often referred to
as `fast simulated annealing'. The acceptance function $f(t)=1/(1+t/\tau)$
introduced above is similar to functions used in such papers;
for instance, it corresponds to $\lambda=1$ in Eq.~(28) of
\citet{gielis-maes-99}, and to $q_A=2$ in Eq.~(5) of
\citet{tsallis-stariolo-96}. None of these authors obtained
rate of convergence results for these schemes however.
\citet[Example~3]{tsallis-stariolo-96} did obtain a convergence rate
for $f(t)=1/(1+t)^2$ and showed that this rate is indeed
faster than for classical simulated annealing; the result however assumes
that $\psi_{\max}$ is known and these authors worked exclusively
on a finite set $\Theta$.

We now return to the algorithm and define, for any $\beta>0$
and $x,y\in\Theta$,
$$
   a_\beta(x,y) = f(\beta(\psi(x)-\psi(y))_+).
$$
One step of the above algorithm is then described by a Markov transition kernel
$K_\beta$ defined as
\begin{equation}
\label{Eq:defK}
K_\beta(x,dy) = a_\beta(x,y) K(x,dy)
  + \left(  1-\int_\Theta a_\beta(x,z)\,K(x,dz)  \right) \delta_x(dy) .
\end{equation}
Thus, assuming that the initial point $\theta_0$ is random and drawn
from some probability distribution $\eta_0$ on $\Theta$,
the sequence $(\theta_n)_{n\geq 0}$ is an inhomogeneous Markov chain with
initial law $\eta_0$ and transition kernels $(K_{\beta_n})_{n\geq 0}$;
more precisely, for any $n$, $\theta_{n+1}$ has conditional distribution
$K_{\beta_n}(\theta_n,\cdot)$. 

We will suppose that $\Theta$ is equipped with its Borel $\sigma$-field  
$\mathcal{B}(\Theta)$, and we will also assume that the Markov transition 
kernel $K$ satisfies the following condition.
\begin{hypothesis}\label{Hyp:mixing}
There exists $\epsilon_K >0$ and a probability measure $\lambda$ on
$(\Theta,\mathcal{B}(\Theta))$ such that
$$
\forall (x,A) \in \Theta \times \mathcal{B}(\Theta):\;
\epsilon_K \lambda(A) \leq  K(x,A)\leq \frac{1}{\epsilon_K} \lambda(A) .
$$
\end{hypothesis}
\noindent
Of course, Hypothesis \ref{Hyp:mixing} is easier to fulfil if
$\Theta$ is compact or bounded.

Regarding the function $\psi$, we also make some assumptions.
Put, for any $\epsilon>0$ and $a<b$,
$$
U^{\epsilon,a,b}=\{x\in \Theta:\psimax-b-\epsilon< \psi(x) \leq
  \psimax-a-\epsilon\} .
$$

\begin{hypothesis}\label{Hyp:osc}
The oscillations of $\psi$ are bounded, that is,
$$
\osc(\psi):=\sup_{x,y \in \Theta}|\psi(x)-\psi(y)|<\infty .
$$
\end{hypothesis}

\begin{hypothesis}\label{Hyp:gradpsi}
Either one of the following two assumptions holds true.
\begin{itemize}
\item[(i)] For all $\epsilon>0$ small enough
  there are numbers $C_0(\epsilon)>0$ and $\epsilon'>0$
  such that for all $\delta>0$,
$$
\lambda(U^{\epsilon,\delta i,\delta(i+1)})\leq C_0 \delta
\quad\mbox{or}\quad i\delta \geq \epsilon'.
$$
\item[(ii)]  The function $\psi$ has a single global maximum,
  $\theta_{\max}$ say, located in the interior of $\Theta$
  (which is thus non-empty). The probability measure $\lambda$ is absolutely
  continuous with respect to Lebesgue measure and its density is locally
  bounded. The function $\psi$ is $C^3$ in
  $\{x:\psi(x)>\psimax - \epsilon''\}$, which is a neighbourhood of
  $\theta_{\max}$ (for some $\epsilon''>0$), and the quadratic form
  $\psi''(\theta_{\max})$ is negative definite.
\end{itemize}
\end{hypothesis}
The attentive reader will notice that one could replace the assumption of
a unique maximum by an assumption that there are a finite number of maxima,
and that one could replace (ii) above by some more sophisticated assumptions
on the derivatives of $\psi$. This requires a higher level of technicality
but the whole proof would contain the same ideas and this is why we write
the assumptions in this way.

\section{Rate of convergence}\label{sec:rate}

Throughout the remainder of the paper we take $\beta_n=n^\alpha\vee 1$
for some  $0<\alpha<1$. The choice of this
particular sequence will be explained in
Remark~\ref{Rem:betachoice}. We denote by $(\theta_n)_{n\geq 0}$ the
sequence produced by the annealing scheme for this cooling schedule.
The main result of the present section is the following.

\begin{theorem}
\label{Th:rate}
Let $M=\osc(\psi)\vee\osc(\psi)^\alpha$ and suppose that
$g(t^\alpha)/t \rightarrow 0$ as $t \rightarrow \infty$.
Then for all $\epsilon>0$ small enough
(if Hypothesis~\ref{Hyp:gradpsi}(i) holds) or
$0<\epsilon\leq\epsilon''$ (if Hypothesis~\ref{Hyp:gradpsi}(ii) holds),
there exists a $C_\epsilon>0$ depending on $\epsilon$
such that for all $n$,
\begin{eqnarray*}
\lefteqn{\p(\psi(\theta_n)\leq \psimax-\epsilon)} \hspace*{10mm} \\
& \leq & C_{\epsilon} \left(  \frac{g(M (n+1)^\alpha)^2}{n}
  + \frac{1}{n^\alpha}
  \left( 1 +\int_0^{1+n^\alpha \osc(\psi)} f(t)\, dt \right)\right)
   +f(\epsilon' n^\alpha) \\
&& \text{under Hypothesis \ref{Hyp:gradpsi}(i)}, or \\
\lefteqn{\hspace*{-10mm}\p(\psi(\theta_n)\leq \psimax-\epsilon)} \\
& \leq & C_{\epsilon} \left(  \frac{g(M(n+1)^\alpha)^2}{n}
  + \frac{1}{n^\alpha} \int_0^{n^\alpha \osc(\psi)^2} f(t)\, dt \right)
   +f((\epsilon''-\epsilon) n^\alpha) \\
&& \text{under Hypothesis \ref{Hyp:gradpsi}(ii)}.
\end{eqnarray*}
\end{theorem}

\begin{corollary}\label{Cor:optimal}
Choosing $\alpha=1/3$ and $g(t)=1+t/\tau$, where $\tau>0$ is arbitrary,
the bounds of Theorem~\ref{Th:rate} are $C_\epsilon C n^{-1/3}(1+\log n)$.
\end{corollary}

\begin{remark}\label{Rem:optimalchoice}
If we want to have terms of the same order in the bounds of
Theorem~\ref{Th:rate}, we see that $g(M n^\alpha)^2/n$ and
$f(\epsilon' n^\alpha)$ (or $f((\epsilon''-\epsilon) n^\alpha)$,
depending on the case)
should be of the same order. Thus $f(t)$ should be of order $t^{-1/3}$ as
$t\to\infty$. With this choice all terms in the bound have the same order,
and so there is something optimal to it. With our inequalities,
it does not seem possible to have a better rate.
\end{remark}

\begin{remark}\label{Rem:freetemp}
In Corollary \ref{Cor:optimal} there is a parameter $\tau>0$
which can be chosen arbitrarily. This parameter plays the role of a
temperature like in classical simulated annealing and can be tuned by
the user to optimise convergence. On the contrary to classical
simulated annealing there is, theoretically, no restriction on $\tau$.
\end{remark}

Before going into the proof of these results, we will proceed through
some technical lemmas. First we however give some additional
notation.
The total variation distance $\|\mu-\nu\|_{\TV}$ between two probability
measures $\mu$ and $\nu$ is defined as $\sup_A|\mu(A)-\nu(A)|$,
where the supremum is taken over the $\sigma$-field on which the
measures are defined. The set of probability measures on
$(\Theta,\mathcal{B}(\Theta))$ will be denoted by $\mathcal{P}(\Theta)$.

\begin{lemma}\label{Lem:mixbeta}
For any $\beta>0$ it holds that
$$
\forall (x,A) \in \Theta \times \mathcal{B}(\Theta):\; K_\beta(x,A)
  \geq \epsilon_K f(\beta \osc(\psi)) \, \lambda(A).
$$
\end{lemma}

\begin{corollary}\label{Cor:contraction}
The preceding lemma and Dobrushin's theorem (see \citealp{dobrushin-56},
or \citealp{delmoral-guionnet-01}) imply that for any $\beta >0$ and
any probability measures $\mu$ and $\nu$ on $\Theta$,
$$
\Vert \mu K_\beta - \nu K_\beta \Vert_{\TV}
  \leq (1-\epsilon_K f(\beta \osc(\psi))) \Vert \mu-\nu\Vert_{\TV}.
$$
\end{corollary}

\begin{proof}[Proof of Lemma \ref{Lem:mixbeta}]
Take $\beta>0$ and $(x,A) \in \Theta \times \mathcal{B}(\Theta)$. Then
\begin{eqnarray*}
K_\beta(x,A) &\geq& \int_A a_\beta(x,y) \, K(x,dy)\\
& \geq & \int_A f(\beta \osc(\psi)) \epsilon_K \, \lambda(dy)\\
& = & \epsilon_K f(\beta \osc(\psi)) \lambda(A).
\end{eqnarray*}
\end{proof}

The above corollary implies that for any $\mu$, the sequence
$(\mu K_\beta^n)_{n\geq 0}$ is a Cauchy sequence in total variation
norm. Thus there exists a total variation limit
\citep[cf.][p.~232]{lindvall-2002}, which we denote by $\mu_\beta$.
This probability measure is invariant for $K_\beta$, and it does not 
depend on the particular choice of the initial distribution $\mu$. 
It is hence the unique invariant distribution of $K_\beta$.

The convergence of simulated annealing hinges on the fact that the law of
$\theta_n$, which we denote by $\eta_n$, is close to $\mu_{\beta_n}$,
and that for
large $\beta_n$ the measure $\mu_{\beta_n}$ is concentrated on the regions
where $\psi$ is large. This concentration is the subject of the next lemma.
We set
$$
U^\epsilon=\{x \in \Theta:\psi(x)>\psimax-\epsilon\},
$$
and $\uec$ is its complement in $\Theta$.

\begin{lemma}\label{Lem:stationary}
For all $\beta>0$ and $\epsilon>0$ small enough
(if Hypothesis~\ref{Hyp:gradpsi}(i) holds) or
$0<\epsilon\leq\epsilon''$ (if Hypothesis~\ref{Hyp:gradpsi}(ii) holds),
there is a constant $C_\epsilon$ depending on $\epsilon$ such that
$$
\begin{array}{rcll}
\mu_\beta(\uec) & \leq &
  \ds \frac{C_\epsilon}{\beta}  \left(
  1+ \int_0^{1+\beta\osc(\psi)} f(t) \,dt\right)
  + f( \beta\epsilon')
  & \text{under Hypothesis \ref{Hyp:gradpsi}(i),}\\[5mm]
\mu_\beta(\uec) & \leq &
  \ds \frac{C_\epsilon}{\beta} \int_0^{\beta\osc(\psi)^2} f(t)\,dt
   + f( \beta(\epsilon''-\epsilon))
   & \text{under Hypothesis \ref{Hyp:gradpsi}(ii)}.
\end{array}
$$
\end{lemma}

\begin{proof}
Fix $\beta>0$ and $\epsilon$ in the appropriate range. We have
\begin{eqnarray*}
\mu_\beta(U^{\epsilon,c}) & = & \mu_\beta K_\beta(U^{\epsilon,c})\\
& = & \iint_{x\in\uec,y\in \uec } \mu_\beta(dx)\,K_\beta(x,dy)\\
&&~~~~   + \iint_{x\in \ue,y \in \uec} \mu_\beta(dx)\,K_\beta(x,dy).
\end{eqnarray*}
For $x\in\uec$ and $y\in\ue$, $K_\beta(x,dy)=K(x,dy)$. Thus
the first integral above can be bounded as
\begin{eqnarray*}
\iint_{x\in\uec,y\in \uec } \mu_\beta(dx)\,K_\beta(x,dy)
  & = & \int_{x\in \uec} \mu_\beta(dx)
        \left( 1-\int_{y \in \ue} K (x,dy)  \right)\\
& \leq & \int_{x\in \uec } \mu_\beta(dx)
   \left( 1 - \int_{y\in \ue} \epsilon_K \,\lambda(dy)  \right)\\
& = & (1-\epsilon_K \lambda(\ue)) \mu_\beta(\uec).
\end{eqnarray*}
Similarly, for the second integral,
\begin{eqnarray*}
\lefteqn{\iint_{x\in\ue,y\in \uec } \mu_\beta(dx)\,K_\beta(x,dy)}
  \hspace*{10mm} \\
& = & \iint_{x\in \ue, y \in \uec} \mu_\beta(dx) a_\beta(x,y) \,K(x,dy) \\
& \leq  & \iint_{x\in \ue,y\in \uec} \mu_\beta(dx) a_ \beta(x,y)
      \frac{1}{\epsilon_K}\,\lambda(dy) \\
& \leq & \frac{1}{\epsilon_K} \int_{y \in \uec}
   f(\beta(\psimax-\epsilon-\psi(y)))\, \lambda(dy),
\end{eqnarray*}
so that
$$
\mu_\beta(\uec) \leq \frac{1}{\epsilon_K^2 \lambda(\ue)}
   \int_{y \in \uec} f(\beta(\psimax-\epsilon-\psi(y))) \, \lambda(dy) .
$$

To finish the proof we will now bound the above integral as in the
statement of the lemma.
If Hypothesis~\ref{Hyp:gradpsi}(i) holds, take $\delta=1/\beta$
and proceed as
\begin{eqnarray*}
\lefteqn{\int_{y \in \uec} f(\beta(\psimax-\epsilon-\psi(y))) \, \lambda(dy) 
 \leq  \sum_{i=0}^{\osc(\psi)/\delta} f(i) \lambda(\uei)}
  \hspace*{20mm} \\
& \leq & \sum_{i=0}^{\osc(\psi)/\delta} C_0(\epsilon) \delta f(i)
   + \sum_{i\geq \epsilon'/\delta}f(i) \lambda(\uei)\\
& \leq & C_0(\epsilon) \delta \left( f(0)
   + \int_0^{1+\osc(\psi)/\delta} f(t) \, dt  \right)
   + f\left(\frac{\epsilon'}{\delta}\right).
\end{eqnarray*}

If Hypothesis~\ref{Hyp:gradpsi}(ii) holds we employ Morse's lemma
\cite[see e.g.][Theorem~4.2.12]{berger-gostiaux-88}
to make a change of variables in
$\{x:\psi(x)>\psimax - \epsilon''\}$ such that with some bounded function
$\xi$ (that only depends on $\psi$),
\begin{eqnarray*}
\lefteqn{\int_{y \in \uec} f(\beta(\psimax-\epsilon-\psi(y))) \, \lambda(dy)}
  \hspace*{20mm} \\
& \leq & \int_{\psi(y)>\psimax-\epsilon'',y\in \uec}
   f(\beta(\psimax-\epsilon-\psi(y))) \, \lambda(dy) \\
&& + \int_{\psi(y)\leq \psimax-\epsilon''} f(\beta(\epsilon''-\epsilon))
    \, \lambda(dy) \\
& \leq & \int_{\sqrt{\epsilon}}^{\osc(\psi)} f(\beta (t^2-\epsilon)) \xi(t)\, dt +  f(\beta(\epsilon''-\epsilon))\\
& \leq &  \frac{\|\xi\|_\infty}{2\sqrt{\epsilon}\beta}
   \int_0^{\beta \osc(\psi)^2} f(u) \,du
   +  f(\beta(\epsilon''-\epsilon)),
\end{eqnarray*}
after a change of variable $u=\beta(t^2-\epsilon)$.
\end{proof}

\begin{remark} \label{Rem:betachoice}
In the following we will show that $\eta_n$ is close to $\mu_{\beta_n}$.
Using Lemma~\ref{Lem:stationary} to bound $\mu_{\beta_n}(\uec)$,
we obtain a bound larger than $1/\beta_n$. We would like to compare
$\eta_n(\uec)$ to a power of $n$, so it is natural at this point to take,
for some $\alpha>0$,
$$
\beta_n=n^\alpha \vee 1.
$$
For technical reasons appearing in the proof of Theorem~\ref{Th:rate},
we need to take $\alpha<1$.
\end{remark}

The law $\eta_n$ approaches $\mu_{\beta_n}$ which becomes increasingly
concentrated on regions where $\psi$ is large,
but at the same time $\beta_n$ is changing.
The following lemma serves us to bound the distance between
$\mu_\beta$ and $\mu_{\beta '}$.

\begin{lemma}\label{Lem:diffbeta}
With $C=(1/\epsilon_K)\sup_{t\geq 0}|tf'(t)|$ it holds that for
any $\beta'>\beta>0$, 
$$
\Vert \mu_\beta - \mu_{\beta'} \Vert_{\TV}
  \leq C g(\beta \osc(\psi)) \left( \frac{\beta'}{\beta} -1 \right) .
$$
\end{lemma}

\begin{proof}
We have, using Corollary \ref{Cor:contraction},
\begin{eqnarray*}
\lefteqn{\Vert \mu_\beta -\mu_{\beta'} \Vert_{\TV} \leq
  \Vert \mu_{\beta}K_{\beta} - \mu_{\beta'}K_{\beta} \Vert_{\TV}
  + \Vert \mu_{\beta'}K_{\beta} - \mu_{\beta'}K_{\beta'} \Vert_{\TV} }
   \hspace*{10mm}\\[3mm]
& \leq & (1-\epsilon_K f(\beta \osc(\psi)))
  \Vert \mu_\beta -\mu_{\beta'} \Vert_{\TV}
  + \sup_{\mu \in \mathcal{P}(\Theta)}
  \Vert \mu K_{\beta} - \mu K_{\beta'} \Vert_{\TV}.
\end{eqnarray*}

Pick $\mu\in\mathcal{P}(\Theta)$. We may construct two coupled samples from
$\mu K_\beta$ and $\mu K_{\beta'}$ respectively by first sampling
$x$ from $\mu$, then sampling $y$ from $K(x,\cdot)$, sampling $U$
from the uniform distribution on $(0,1)$ and finally accepting the
proposal $y$ if $U\leq \alpha_\beta(x,y)$ or $U\leq \alpha_{\beta'}(x,y)$
respectively. Similarly to Appendix~\ref{app:coupling} we may then
conclude that
\begin{eqnarray*}
\Vert \mu K_{\beta} - \mu K_{\beta'} \Vert_{\TV}
  & \leq & \sup_{x,y \in \Theta} |a_{\beta}(x,y)-a_{\beta'}(x,y)|\\
  & \leq &  \sup_{0\leq u \leq \osc(\psi)}
      \left|f(\beta u ) - f(\beta' u)\right|\\
  & \leq & \sup_{0\leq u \leq \osc(\psi)}
    ( |f'(\xi u)| \beta u )
    \left(\frac{\beta'}{\beta} -1 \right),
\end{eqnarray*}
where $\xi=\xi(u)$ is a point between $\beta$ and $\beta'$. Since
$f$ is assumed convex and non-increasing, and hence $|f'|$ non-increasing,
it holds that
$|f'(\xi u)|\leq |f'(\beta u)|$. We thus arrive at the bound
$$
\Vert \mu_\beta - \mu_{\beta'} \Vert_{\TV}
\leq \frac{1}{\epsilon_K f(\beta \osc(\psi))}
  \sup_{0\leq u \leq \osc(\psi)} (|f'(\beta u)| \beta u)
  \left(\frac{\beta'}{\beta} -1 \right) .
$$
Since $|tf'(t)|$ is assumed bounded, the proof is complete.
\end{proof}

\begin{proof}[Proof of Theorem \ref{Th:rate}]
Set $\Delta_n = \Vert \eta_n - \mu_{\beta_n} \Vert_{\TV}$. If we can
prove the inequality
\begin{equation}\label{eq:Inineq}
\Delta_n \leq C \frac{g(M (n+1)^\alpha)^2}{n},
\end{equation}
the result will follow from Lemma~\ref{Lem:stationary} and,
in case of Hypothesis~\ref{Hyp:gradpsi}(i), the bound
$g\geq 1$.

In order to prove (\ref{eq:Inineq}) the assumption
$g(x^\alpha)/x\to 0$ as $x\to\infty$ will be instrumental.
We start by deriving a recursive bound for $\Delta_n$.
By Corollary~\ref{Cor:contraction} and Lemma~\ref{Lem:diffbeta} we have,
for all $n$,
\begin{eqnarray*}
\Delta_{n+1} & \leq
  & \Vert \eta_n K_{\beta_n} - \mu_{\beta_n}K_{\beta_n} \Vert_{\TV}
  + \Vert \mu_{\beta_n} - \mu_{\beta_{n+1}} \Vert_{\TV}\\[2mm]
& \leq & (1-\epsilon_K f(\beta_n \osc(\psi)))
  \Vert\eta_n - \mu_{\beta_n} \Vert_{\TV}
 + C\left(\frac{\beta_{n+1}}{\beta_n} -1 \right) g(n^\alpha \osc(\psi))
  \\[2mm]
& \leq & (1-\epsilon_K f(\beta_n \osc(\psi))) \Delta_n 
 + C \frac{g(n^\alpha \osc(\psi))}{n+1} .
\end{eqnarray*}
Iterating this recursion yields
\begin{eqnarray*}
\Delta_{n+1} & \leq  & 
  \sum_{q=1}^{n} \prod_{k=q+1}^n (1-\epsilon_K f(\beta_k \osc(\psi)))
  \times C \frac{g(q^\alpha \osc(\psi))}{q+1} \\
& + & \prod_{k=1}^n (1-\epsilon_K f(\beta_k \osc(\psi)))
 \times \|\eta_1 - \mu_{\beta_1}\|_{\TV},
\end{eqnarray*}
where an empty product (when $q=n$) is interpreted as unity.

Define $F$ such that $F'(x)=f(x^\alpha)$. Then for $1\leq q \leq n-1$,
\begin{eqnarray*}
\lefteqn{\log \prod_{k=q+1}^n (1-\epsilon_K f(\beta_k \osc(\psi)))
 \leq -\sum_{k=q+1}^{n} \epsilon_K f(\osc(\psi) k^\alpha)}
  \hspace*{30mm} \\
& \leq & -\epsilon_K \int_{q+1}^{n+1} f(\osc(\psi) x^\alpha) \,dx \\
& = & - \frac{\epsilon_K}{\osc(\psi)}
  (F(\osc(\psi)(n+1))-F(\osc(\psi)(q+1))) .
\end{eqnarray*}
For $q=n$ this is an equality.
Putting $C_1=\epsilon_K/\osc(\psi)$ and  $C_2= \osc(\psi)$ we thus obtain
\begin{eqnarray*}
\Delta_{n+1} & \leq 
& e^{- C_1 F(C_2(n+1))} \sum_{q=1}^{n} e^{C_1 F(C_2(q+1))}
  C \frac{g(C_2 q^\alpha)}{q+1}\\
& + & 2 e^{- C_1 (F(C_2(n+1))-F(C_2))} \\
& \leq & C e^{-C_1 F(C_2(n+1))} \int_1^{n+1}e^{C_1F(C_2(x+1))}
  \frac{g(C_2 x^\alpha)}{x} \, dx \\
& + & 2 e^{- C_1 (F(C_2(n+1))-F(C_2))} .
\end{eqnarray*}

Denote the integral on the right-hand side by $I_{n+1}$. First we notice
that since $g\geq 1$, $I_{n+1}\to\infty$ as $n\to\infty$. Next we rewrite
this integral as
$$
 I_{n+1}=\int_1^{n+1}e^{C_1F(C_2(x+1))} C_1C_2 f(C_2^\alpha (x+1)^\alpha)
 \times \frac{g(C_2^\alpha (x+1)^\alpha)}{C_1C_2}
    \frac{g(C_2 x^\alpha)}{x} \, dx,
$$
where $C_1C_2f(C_2^\alpha(x+1)^\alpha)$ is the derivate
of the exponent.
By partial integration, integrating the first factor of the integrand
above and dropping all negative contributions (recall that $g$ and $g'$
are non-negative), we obtain the bound
\begin{eqnarray}
I_{n+1}
& \leq & \left[ e^{C_1F(C_2(x+1))} \frac{g(C_2^\alpha (x+1)^\alpha)}{C_1C_2}
    \frac{g(C_2 x^\alpha)}{x} \right]_1^{n+1}  \nonumber \\
& + & \int_1^{n+1} e^{C_1F(C_2(x+1))} \frac{g(C_2^\alpha (x+1)^\alpha)}{C_1C_2}
    \frac{g(C_2 x^\alpha)}{x^2} \,dx. \label{Eq:remainderint}
\end{eqnarray}
Denote the integral on the right-hand side of (\ref{Eq:remainderint})
by $I'_{n+1}$. This integral is similar to $I_{n+1}$, the difference
being that the integrand is multiplied by a constant times
$g(C_2^\alpha (x+1)^\alpha)/x$. Since this ratio tends to
zero as $x\to\infty$, and since $I_{n+1}\to\infty$ (as noted above),
it holds that for any $0<\kappa<1$,
$I'_{n+1}\leq\kappa I_{n+1}$ for sufficiently large $n$. Hence
\begin{eqnarray*}
I_{n+1} & \leq &
\frac{1}{1-\kappa}
  \left[ e^{C_1F(C_2(x+1))} \frac{g(C_2^\alpha (x+1)^\alpha)}{C_1C_2}
    \frac{g(C_2 x^\alpha)}{x} \right]_1^{n+1} \\
& \leq & C e^{C_1F(C_2(n+2))} \frac{g(M (n+2)^\alpha)^2}{n+1}
\end{eqnarray*}
for sufficiently large $n$;
recall that $M=C_2\vee C_2^\alpha$. Summing up thus far, we have shown that
\begin{eqnarray}
\Delta_{n+1} & \leq & C \left(
  e^{C_1 (F(C_2(n+2))-F(C_2(n+1)))} \frac{g(M (n+2)^\alpha)^2}{n+1}
  + e^{-C_1 F(C_2(n+1))} \right) \nonumber \\
& \leq & C \left(  \frac{g(M (n+2)^\alpha)^2}{n+1}
  + e^{-C_1 F(C_2(n+1))} \right), \label{Eq:Deltatwotermbound}
\end{eqnarray}
where the second inequality follows as $F'$ is bounded.

Now take an arbitrary $m>0$. Since $g(x^\alpha)/x\to 0$ there is an
$x_m>0$ such that $g(x^\alpha)/x\leq 1/m$ for $x\geq x_m$, or,
equivalently, $f(x^\alpha)\geq m/x$ for $x\geq x_m$. Integrating
this inequality yields $F(x)-F(x_m)\geq m\log(x/x_m)$, so that
$$
e^{-C_1 F(x) + C_1 F(x_m)} \leq \left( \frac{x}{x_m}\right)^{-mC_1}
$$
for $x\geq x_m$.
Picking $m$ such that $mC_1=1$ we see that as $n\to\infty$, the
second term on the right-hand side of (\ref{Eq:Deltatwotermbound})
is of smaller order than the first one. We conclude that
$$
\Delta_{n+1} \leq C \frac{g(M (n+2)^\alpha)^2}{n+1},
$$
which is (\ref{eq:Inineq}).
\end{proof}


\section{Simulated annealing on a function that cannot
  be computed exactly} \label{Sec:unfixed-observations}

In this section we assume that the function $\psi$ to be maximised cannot
be computed explicitly, but that we have available an approximation
to it. This approximation, denoted by $\psi^N$, can be stochastic,
based on Monte Carlo
procedures; the next section shows such an example.
The precision of the approximation, stochastic or not, is indexed
by an integer-valued parameter $N$, and the larger the $N$, the
better the approximation. This parameter can be, for instance,
the number of replications in a Monte Carlo method. The following
hypothesis makes precise the quality of the approximation.


\begin{hypothesis}\label{Hyp:approx}
For all $N\geq 1$ we can compute a deterministic or stochastic
approximation $\psi^N$ of $\psi$ such that
$$
\E|\psi^N(x)-\psi(x)|\leq \frac{a_1}{\sqrt{N}}
\quad\mbox{for all $x\in\Theta$},
$$
and, almost surely,
$$
|\psi^N(x)-\psi^N(y)|\leq 2 \osc(\psi)
\quad\mbox{for all $x,y\in\Theta$}.
$$
\end{hypothesis}

We suppose that this hypothesis holds true in all of the following.
The attentive reader will notice that the second of the above assumptions
can be replaced by the existence of a constant $C$ such that,
almost surely, $|\psi^N(x)-\psi^N(y)|\leq C$ for all $x,y\in\Theta$.
In the case of approximation by a sample mean of i.i.d.\ summands,
the first part of the hypothesis follows from the
Marcinkiewicz-Zygmund inequality;
see Appendix~\ref{app:coupling} for more details.

The sequence $(\beta_n)$, the cooling schedule, is again chosen as
\begin{equation}\label{eq:unfixedbeta}
  \beta_n=n^\alpha \vee 1,
\end{equation}
although below we argue for the choice $\alpha=1/4$ rather than
$\alpha=1/3$ as in the previous section.
We will let the parameter $N$ depend on the iteration
number $n$ as well, $N=N_n$, and we will assume that the increase
is affine in $n$, meaning that $N_n=\lceil N_0 + N_1 n \rceil$
for some numbers $N_0\geq 0$ and $N_1>0$ where $\lceil x \rceil$
denotes rounding $x$ upwards to the nearest integer.
We comment on other choices of $(N_n)_{n\geq 0}$
following the proof of Theorem~\ref{Th:unfixed-observations} below.

We now formalise the simulated annealing procedure in this
modified context. The procedure is again described as a random
sequence, denoted by $(\bt_n)_{n\geq 0}$, with $\bt_0$ sampled
from the law $\eta_0$ (as is $\theta_0$). The function $g$ is
chosen as in Corollary~\ref{Cor:optimal}, and $(\bt_n)$
evolves as follows.

\begin{itemize}
\item[(c1)] In stage $n$, given the current state $\bt_n$,
  sample a new proposed position $Z$ from $K(\bt_n,\cdot)$.
\item[(c2)] Set $\bt_{n+1}=Z$ with probability
  $$
      f(\beta_n(\psi^{N_n}(\bt_n)-\psi^{N_n}(Z))_+)
  $$
  and $\bt_{n+1}=\bt_n$ otherwise.
\end{itemize}
This procedure requires some comments. In step~(c2), $\psi$ is approximated
at two points, $\theta_n$ and $Z$. In the case of random approximations
it is unimportant whether these two evaluations are
independent or not, as we shall see below, but it
is important that they are independent of approximations computed
in previous steps (smaller $n$) of the algorithm. The reason for this
is that, if such independence holds, the sequence $(\bt_n)$ forms
a Markov chain, and this Markov chain is the object of our study.
Moreover, the additional randomness in step~(c2) associated with
the phrases `sample a new proposed position\ldots' and
`with probability\ldots', typically obtained by drawing
random numbers uniformly in $(0,1)$, must be based on two
mutually independent sequences of independent random numbers, also
independent of the function approximations $\psi^N$; this is just as
in the previous annealing schemes however.

In cases where the random function approximations $\psi^N$ are such
that they depend on random variables that are drawn once and for all
and then stay fixed over $n$ (sometimes called `fixed randomness'),
so that $\psi^N$ is fixed at each point in $\Theta$,
we can, as long as $N$ stays fixed too, apply the results of
the previous section to the function $\psi^N$ provided that it satisfies
the regularity assumptions made there. Main questions are then rather whether
these assumptions indeed are satisfied for $\psi^N$, and how well the maximum
of $\psi^N$ and its location approximate those of $\psi$.

We now return to the sequence $(\bt_n)$. As noted above, this
sequence is an (inhomogeneous) Markov chain.
For any $\beta>0$ and $N\geq 1$, we define the function 
 $$
 a^N_\beta(x,y) = f(\beta(\psi^N(x)-\psi^N(y))_+) .
$$
For fixed $x$ and $y$ this is indeed a random variable, the randomness
coming from the evaluations $\psi^N(x)$ and $\psi^N(y)$. We write
$\Ep$ for the expectation with respect to the random variables used to
compute $\psi^N$ at a point for some approximation index $N$,
and $\pp$ for the corresponding probability. The kernels
$K^{N_n}_{\beta_n}$ of $(\bt)_{n\geq 0}$, defined by 
$$
K^{N_n}_{\beta_n}(x,A) =  \p(\bt_{n+1}\in A\mid\bt_n=x)  
$$
for any $x\in\Theta$ and $A\in \mathcal{B}(\Theta)$, can then be expressed as
\begin{equation}\label{Eq:defKN}
  K^{N}_{\beta}(x,dy) = \Ep\!\! \left[ a^{N}_{\beta}(x,y) K(x,dy)
    + \left(  1-\int_\Theta a^{N}_{\beta}(x,z)\,K(x,dz)  \right)
    \delta_x(dy)\right].
\end{equation}

The final assumption we make before stating the main result of this
section is the following.

\begin{hypothesis}\label{Hyp:bias}
There is a constant $C_K$ such that for all $\beta>0$ and $N'>N\geq 1$,
\begin{equation}\label{eq:biasbound}
\sup_{\mu \in \mathcal{P}(\Theta)}
  \| \mu K_\beta^N-\mu K_\beta^{N'}\|_{\TV}
  \leq C_K \beta\frac{N'-N}{N}.
\end{equation}
\end{hypothesis}
In Appendix~\ref{app:coupling} we discuss this condition in detail
for approximations obtained as sample means of i.i.d.\ random variables,
and for approximations obtained using so-called particle filters.
It turns out that Hypothesis~\ref{Hyp:bias} can often be verified
through a coupling argument; that is, we couple the approximations
$\psi^N$ and $\psi^{N'}$ in a suitable way. We notice that
by such an argument it also follows that provided Hypothesis~\ref{Hyp:approx}
holds, one can bound the left-hand side of (\ref{eq:biasbound}) by
a constant times $\beta/\sqrt{N}$; the actual assumption above is
thus stronger.

\begin{theorem}\label{Th:unfixed-observations}
Assume that $\beta_n$ is as in (\ref{eq:unfixedbeta}) with
$\alpha<1/2$ and that $N_n$ increases linearly with $n$.
Then under Hypotheses \ref{Hyp:mixing}, \ref{Hyp:osc}, \ref{Hyp:gradpsi},
\ref{Hyp:approx} and \ref{Hyp:bias},
for all $\epsilon>0$ small enough
(if Hypothesis~\ref{Hyp:gradpsi}(i) holds) or
$0<\epsilon\leq\epsilon''$ (if Hypothesis~\ref{Hyp:gradpsi}(ii) holds),
there exists a constant
$C'_\epsilon$ depending on $\epsilon$ such
that
\begin{equation}\label{Eq:boundgeneral}
\p(\psi(\bt_n ) \leq \psimax - \epsilon)
  \leq C_\epsilon' (n^{-\alpha}  \log n \vee n^{3\alpha-1})
\end{equation}
for sufficiently large $n$.
\end{theorem}
Equating the two powers of this bounds leads to
$\alpha=1/4$ as the optimal choice, with corresponding rate of convergence
$n^{-1/4}\log n$.

The proof of Theorem~\ref{Th:unfixed-observations} 
is very similar to the proof of
Theorem~\ref{Th:rate}; before going into it, we will proceed through
some technical lemmas. The following results can be shown exactly
in the same manner as Lemma \ref{Lem:mixbeta} and
Corollary~\ref{Cor:contraction}.

\begin{lemma} \label{Lem:mixbetageneral}
For any $\beta>0$ and $N\geq 1$ it holds that
$$
\forall (x,A) \in \Theta \times \mathcal{B}(\Theta):\;
K^N_\beta(x,A) \geq \epsilon_K\lambda(A) f(2\beta \osc(\psi)) .
$$
\end{lemma}

\begin{corollary}\label{Cor:contractiongeneral}
For all $\beta>0$, $N\geq 1$ and any probability measures $\mu$ and $\nu$
on $\Theta$,
$$
\Vert \mu K^N_\beta - \nu K^N_\beta \Vert_{\TV}
  \leq (1- \epsilon_K f(2\beta \osc(\psi)))\Vert \mu-\nu\Vert_{\TV}.
$$
\end{corollary}
We point out, in particular, that these results hold true regardless
of whether the two function approximations required for computing
$a_\beta^N(x,y)$ are independent or not.

The results imply that for any $\beta>0$ and $N\geq 1$,
the kernel $K^N_\beta$ has a unique
stationary distribution, which we denote by $\mu^N_\beta$. We will show that
under certain conditions, $\mu^N_{\beta}$ is concentrated around the
maximum of $\psi$.

\begin{lemma}\label{Lem:devmeasure}
For all $\beta>0$, $\epsilon>0$ small enough
(if Hypothesis~\ref{Hyp:gradpsi}(i) holds) or
$0<\epsilon<\epsilon''$ (if Hypothesis~\ref{Hyp:gradpsi}(ii) holds)
and $N\geq 1$ such that $N\geq\beta^2$ and
$N\geq(8a_1/\epsilon )^2$, there is a constant $C''_\epsilon$ depending
on $\epsilon$ but not on $N$ such that
$$
\mu_{\beta}^N(\uec) \leq C''_\epsilon \frac{1+\log\beta}{\beta} .
$$
\end{lemma}

\begin{proof}
We proceed as in Lemma \ref{Lem:stationary} and thus write
\begin{eqnarray*}
\lefteqn{\mu^N_{\beta}(U^{\epsilon,c}) =
  \mu^N_{\beta} K^N_{\beta}(U^{\epsilon,c})} \hspace*{5mm} \\[2mm]
& = & \iint_{x\in\uec,y\in \uec } \mu^N_{\beta}(dx)\,K^N_{\beta} (x,dy) 
+ \iint_{x\in \ue,y \in \uec} \mu^N_{\beta} (dx) \, K^N_{\beta}(x,dy)
\end{eqnarray*}
and
$$
\iint_{x\in\uec,y\in \uec } \mu^N_{\beta}(dx)\, K^N_{\beta}(x,dy)
 = \int_{x\in \uec} \mu^N_{\beta}(dx) \left( 1-\int_{y \in \ue} K^N_{\beta} (x,dy)  \right).
$$
For $x\in\uec$ it holds that
\begin{eqnarray*}
\lefteqn{\int_{y \in \ue} K^N_{\beta} (x,dy) 
=  \Ep \left( \int_{y\in\ue} f(\beta(\psi^N(x)-\psi^N(y))_+)\,
   K(x,dy)  \right) }\hspace*{10mm} \\
& \geq & \int_{y\in \ue} \epsilon_K \pp(\psi^N(y)\geq \psi^N(x)) \,\lambda(dy)\\
& \geq & \int_{y\in U^{\epsilon/2}} \epsilon_K
  [1-\pp(\psi^N(x)-\psi(x)-\psi^N(y)+\psi(y)\geq \epsilon/2)] \, \lambda(dy)\\
& \geq  & \int_{y\in U^{\epsilon/2}} \epsilon_K
   \left( 1-\frac{4 a_1}{\epsilon \sqrt{N}} \right) \lambda(dy)
   \geq \frac{\epsilon_K}{2} \lambda(U^{\epsilon/2}),
\end{eqnarray*}
where $a_1$ is in Hypothesis~\ref{Hyp:approx} and we used
Markov's inequality and the assumption $N\geq (8a_1/\epsilon)^2$. Hence
$$
\int_{x\in\uec,y\in \uec } \mu^N_{\beta}(dx)\,K^N_{\beta}(x,dy)
  \leq  \mu^N_{\beta}(\uec)\left( 1- \frac{\epsilon_K}{2}
    \lambda(U^{\epsilon/2}) \right)
$$
and
$$
\mu^N_{\beta}(U^{\epsilon,c}) \leq
   \frac{2}{\epsilon_K \lambda(U^{\epsilon/2})}
   \iint_{x\in \ue,y \in \uec} \mu^N_{\beta} (dx) \, K^N_{\beta}(x,dy).
$$

The integral in this bound equals
\begin{eqnarray*}
\lefteqn{\Ep \left( \iint_{x\in\ue,y\in \uec } \mu^N_{\beta}(dx)
   f(\beta(\psi^N(x)-\psi^N(y))_+) \, K(x,dy) \right)}
  \hspace*{3mm} \\
& \leq & \Ep \left(\iint_{x\in\ue,y\in \uec }  \mu^N_{\beta}(dx)
  \frac{1}{\epsilon_K}  f(\beta(\psi^N(x)-\psi^N(y))_+) \, \lambda(dy)  \right)\\
& = & \frac{1}{\epsilon_K} \iint_{x\in\ue,y\in \uec }  \mu^N_{\beta}(dx)
  f(\beta(\psi(x)-\psi(y))_+)\\
&& \hspace*{15mm}\times \Ep \left( \frac{f(\beta(\psi(x)-\psi(y))_+
   + \beta R^N(x,y))}{f(\beta(\psi(x)-\psi(y)))_+} \right) \lambda(dy),
\end{eqnarray*}
where $R^N(x,y)=(\psi^N(x)-\psi^N(y))_+-(\psi(x)-\psi(y))_+$.
Notice that the expression inside the final expectation is bounded
by $g(\beta(\psi(x)-\psi(y))_+)$. Thus the expectation itself,
inserting $g(t)=1+t/\tau$, may be bounded as
\begin{eqnarray*}
\lefteqn{\Ep \left( \frac{1 + \frac{\beta}{\tau}(\psi(x)-\psi(y))_+}
  {1 +\frac{\beta}{\tau}(\psi(x)-\psi(y))_+
    + \frac{\beta}{\tau} R^N(x,y)}\right)} \hspace*{5mm}  \\
& \leq & 2\, \pp\left( \left| \frac{\beta}{\tau} R^N(x,y) \right|
  \leq \frac{1+\frac{\beta}{\tau}(\psi(x)-\psi(y))_+)}{2} \right)\\
& + & \left( 1+\frac{\beta}{\tau}(\psi(x)-\psi(y))_+\right)
   \pp\left( \left| \frac{\beta}{\tau} R^N(x,y) \right|
     \geq \frac{1+\frac{\beta}{\tau}(\psi(x)-\psi(y))_+)}{2} \right)\\
& \leq & 2 + \frac{2\beta}{\tau} \Ep|R^N(x,y)|.
\end{eqnarray*}


Now notice that in the expression for $R^N(x,y)$,
$\psi(x)-\psi(y)\geq 0$ for those $x$ and $y$ appearing in the integral.
It is easy to check that for any real $a$ and $b$ such that
$b\geq 0$, $|a_+-b|\leq |a-b|$. Hence
$|R^N(x,y)|\leq |\psi^N(x)-\psi(x)|+|\psi^N(y)-\psi(y)|$,
and the above expectation is thus bounded by
$2+4\beta a_1/(\tau\sqrt{N}) \leq 2+4a_1/\tau$, where we used the
assumption $N\geq \beta^2$.

As in the proof of Lemma~\ref{Lem:stationary} we may conclude that
\begin{eqnarray*}
\lefteqn{\iint_{x\in\ue,y\in \uec } \mu^N_{\beta}(dx)\,K^N_{\beta}(x,dy)}
  \hspace*{10mm} \\
& \leq & \frac{2(1+2a_1/\tau)}{\epsilon_K}
  \int_{y\in\uec} f(\beta(\psi_{\max}-\epsilon-\psi(y))) \, \lambda(dy).
\end{eqnarray*}

Summing up thus far, we have proved that
$$
\mu^N_{\beta}(U^{\epsilon,c})
\leq \frac{4(1+2a_1/\tau)}{\epsilon_K^2 \lambda(U^{\epsilon/2})}
  \int_{y\in\uec} f(\beta(\psi_{\max}-\epsilon-\psi(y))) \, \lambda(dy).
$$
This integral can be bounded just as in the proof of
Lemma~\ref{Lem:stationary}, and with $f(t)=1/(1+t/\tau)$ these bounds
are of order $C_\epsilon''(1+\log\beta)/\beta$.
\end{proof}

We now formulate an analogue of Lemma~\ref{Lem:diffbeta}.
\begin{lemma}\label{Lem:diffbetageneral}
For any $\beta'>\beta>0$,
$$
\Vert \mu^N_\beta - \mu^N_{\beta'} \Vert_{\TV}
  \leq \frac{1+2\beta \osc(\psi)/\tau}{4 \epsilon_K}
  \left( \frac{\beta'}{\beta}-1\right) .
$$
\end{lemma}

\begin{proof}
We have, using Corollary~\ref{Cor:contractiongeneral},
\begin{eqnarray*}
\lefteqn{\Vert \mu^N_\beta - \mu^N_{\beta'} \Vert_{\TV}
  \leq \Vert\mu^N_{\beta} K^N_{\beta} - \mu^N_{\beta'}K^N_{\beta} \Vert_{\TV}
  +  \Vert\mu^N_{\beta'}K^N_{\beta} - \mu^N_{\beta'}K^N_{\beta'} \Vert_{\TV} }
 \hspace*{10mm}\\[3mm]
& \leq & (1- \epsilon_K f(2\beta \osc(\psi)))
  \Vert\mu^N_{\beta} - \mu^N_{\beta'} \Vert_{\TV}
  + \sup_{\mu \in \mathcal{P}(\Theta)}
   \Vert \mu K^N_{\beta} - \mu K^N_{\beta'}\Vert_{\TV} .
\end{eqnarray*}
Using (\ref{Eq:defKN}), Hypothesis~\ref{Hyp:approx} and an argument as
in the proof of Lemma~\ref{Lem:diffbeta}, we find that for all
$\mu\in\mathcal{P}(\Theta)$,
\begin{eqnarray*}
\Vert \mu K^N_{\beta} - \mu K^N_{\beta'}\Vert_{\TV}
& \leq & \sup_{x,y \in \Theta} \Ep |a^N_\beta(x,y)-a^N_{\beta'}(x,y)| \\
& \leq &  \sup_{0\leq u \leq 2 \osc(\psi)} |f(\beta u)-f(\beta' u)|\\
& \leq & \sup_{0\leq u \leq 2 \osc(\psi)} (|f'(\beta u)| \beta u)
   \left(\frac{\beta'}{\beta} -1 \right).
\end{eqnarray*}
With $g(t)=1+t/\tau$, the above supremum is bounded by $1/4$. Thus
$$
\Vert \mu^N_\beta - \mu^N_{\beta'} \Vert_{\TV} \leq
\frac{1}{\epsilon_K f(2\beta \osc(\psi))}
   \frac{1}{4}\left(\frac{\beta'}{\beta} -1 \right)
= \frac{1+2\beta \osc(\psi)/T}{4 \epsilon_K}
  \left( \frac{\beta'}{\beta}-1\right) .
$$
\end{proof}

\begin{proof}[Proof of Theorem \ref{Th:unfixed-observations}]
First we notice that given the assumptions, including $\alpha<1/2$,
Lemma~\ref{Lem:devmeasure} shows that
$\mu^{N_n}_{\beta_n}(\uec)$ is bounded by
$C_\epsilon'' n^{-\alpha}(1+\alpha\log n)$ for sufficiently large $n$.
This term is the first part of the maximum in (\ref{Eq:boundgeneral}).

Next we denote by $\bar\eta_n$ the law of $\bt_n$ and put
$\Delta_n = \Vert \bar\eta_n - \mu^{N_n}_{\beta_n} \Vert_{\TV}$.
Write
\begin{eqnarray}
\Delta_{n+1}  & \leq & 
  \Vert \bar\eta_n K^{N_n}_{\beta_n}
  - \mu^{N_n}_{\beta_n} K^{N_n}_{\beta_n} \Vert_{\TV}
+ \Vert \mu^{N_n}_{\beta_n} - \mu^{N_n}_{\beta_{n+1}} \Vert_{\TV}
+ \Vert \mu^{N_n}_{\beta_{n+1}} - \mu^{N_{n+1}}_{\beta_{n+1}} \Vert_{\TV}
  \nonumber \\[2mm]
& \leq  & (1-\epsilon_K f(2\beta_n \osc(\psi))) \Delta_n
 +\frac{1+2\beta_n \osc(\psi)/\tau}{4 \epsilon_K}
 \left( \frac{\beta_{n+1}}{\beta_n}-1\right) \nonumber \\
 && + \Vert \mu^{N_n}_{\beta_{n+1}} - \mu^{N_{n+1}}_{\beta_{n+1}} \Vert_{\TV},
 \label{eq:Deltaboundunfixed}
\end{eqnarray}
where we used Corollary~\ref{Cor:contractiongeneral} and
Lemma~\ref{Lem:diffbetageneral} to bound the first two terms.
With our choice of $\beta_n$, the second term on the right-hand side
is of order $n^{\alpha-1}$.

To bound the third term we proceed as in the proof of
Lemma~\ref{Lem:diffbetageneral}; use
Corollary~\ref{Cor:contractiongeneral} to see that for any $\beta$,
$N$ and $N'$,
\begin{eqnarray*}
\lefteqn{\| \mu_\beta^N - \mu_\beta^{N'} \|_{\TV}
\leq
\| \mu_\beta^N K_\beta^N - \mu_\beta^{N'} K_\beta^N\|_{\TV}
+ \| \mu_\beta^{N'} K_\beta^N - \mu_\beta^{N'} K_\beta^{N'} \|_{\TV} }
  \hspace*{10mm} \\[2mm]
& \leq &
(1-\varepsilon_K f(2\beta\mathrm{osc}(\psi)))
   \| \mu_\beta^N - \mu_\beta^{N'}\|_{\TV}
+ \sup_{\mu\in\mathcal{P}(\Theta)}
 \| \mu K_\beta^N - \mu K_\beta^{N'} \|_{\TV}
\end{eqnarray*}
to arrive at
$$
\| \mu_\beta^N - \mu_\beta^{N'} \|_{\TV}
\leq \frac{1}{\varepsilon_K f(2\beta\mathrm{osc}(\psi)))}
\sup_{\mu\in\mathcal{P}(\Theta)}
 \| \mu K_\beta^N - \mu K_\beta^{N'} \|_{\TV}.
$$
Apply this bound with $\beta=\beta_{n+1}$, $N=N_n$ and $N'=N_{n+1}$
to see that the final term of (\ref{eq:Deltaboundunfixed}) is bounded
by a constant times $\beta_{n+1}^2(N_{n+1}-N_n)/N_n$
under Hypothesis~\ref{Hyp:bias}; this ratio is of order 
$n^{2\alpha-1}$ given that $N_n$ is assumed to be affine
in $n$.

Summing up thus far, we have proved that
$$
\Delta_{n+1}  \leq 
(1-\epsilon_K f(2\beta_n \osc(\psi))) \Delta_n
   + \frac{C}{n^{1-2\alpha}}
$$
for some constant $C$. Using this inequality we can show as in the proof of
Theorem~\ref{Th:rate} that for all $\epsilon >0$,
$\Delta_n\leq C'_\epsilon n^{3\alpha-1}$ for some constant
$C'_\epsilon$ depending on $\epsilon$. Indeed,
in the proof of Theorem~\ref{Th:rate}, replace the factor $x$
in the denominator of the expression that forms the integrand
in $I_{n+1}$ by $x^{1-\alpha}$ and proceed from there.
The term $C'_\epsilon n^{3\alpha-1}$ is the second part 
of the maximum in (\ref{Eq:boundgeneral}).
\end{proof}

One may consider other ways of increasing $N_n$, for instance
$N_n = \lceil N_0+N_1n^\delta \rceil$ for some $\delta>0$. 
For $\delta>1$ the expression $(N_{n+1}-N_n)/N_n$ is then still
of order $n^{-1}$ however, so there is no improvement in the
proof of Theorem~\ref{Th:unfixed-observations} compared to
the case of affine increase. For $\delta<1$ the above expression
is of order $n^{-\delta}$, since the $N_n$ are integer-valued.
The bound corresponding to (\ref{Eq:boundgeneral}) then becomes
of order $n^{-\alpha}  \log n \vee n^{3\alpha-\delta}$, with the
optimal $\alpha$ being $\delta/4$.

The above seems to suggest that the rate $n^{-1/3}\log n$ of
Section~\ref{sec:rate} is unobtainable when the function $\psi$ is
approximated. This is not the case however, but it requires
a slightly different approach to analysis than above, 
and also typically a faster increase of $N_n$. In the proof of
Theorem~\ref{Th:unfixed-observations} we compared
$\bar\eta_n$ to $\mu^{N_n}_{\beta_n}$. Consider instead comparing
to $\mu_{\beta_n}$, as in the proof of Theorem~\ref{Th:rate},
and write
\begin{eqnarray*}
\bar\eta_{n+1} - \mu_{\beta_{n+1}} & = &
\bar\eta_n K_{\beta_n}^{N_n} - \mu_{\beta_n}K_{\beta_n}^{N_n} \\
& + & \mu_{\beta_n}K_{\beta_n}^{N_n} - \mu_{\beta_{n+1}}K_{\beta_n}^{N_n}\\
& + & \mu_{\beta_{n+1}}K_{\beta_n}^{N_n} 
     - \mu_{\beta_{n+1}}K_{\beta_{n+1}}^{N_n} \\
& + & \mu_{\beta_{n+1}}K_{\beta_{n+1}}^{N_n} 
    - \mu_{\beta_{n+1}}K_{\beta_{n+1}}.
\end{eqnarray*}
On the right-hand side the  total variation norm of the first difference
is bounded by
$(1- \epsilon_K f(2\beta_n \osc(\psi)))\| \bar\eta_n-\mu_{\beta_n}\|_{\TV}$
(Corollary~\ref{Cor:contractiongeneral}), and the norms of the
remaining differences are bounded by terms of order
$n^{\alpha-1}$ (Lemma~\ref{Lem:diffbeta}), 
$n^{-1}$ (use part of the proof of Lemma~\ref{Lem:diffbetageneral})
and $n^\alpha/N_n^{1/2}$ respectively.
To obtain the order $n^\alpha/N_n^{1/2}$ of the final term
we can couple the kernels $K_\beta$ and $K_\beta^N$
in a way similar to that used in the first part of 
Appendix~\ref{app:coupling}, thus obtaining a bound on the total
variation distance of order 
$\beta\sup_{x\in\Theta}\E_N|\psi_N(x)-\psi(x)|$; by 
Hypothesis~\ref{Hyp:approx} this expression is of order $\beta/N^{1/2}$.
Thus we do not require Hypothesis~\ref{Hyp:bias} for this analysis.

We can now put $\Delta_n=\| \bar\eta_n - \mu_{\beta_n}\|_{\TV}$
and mimic the proof of Theorem~\ref{Th:unfixed-observations}.
To obtain the rate of convergence $n^{-\alpha}\log n$, the norms
of all differences on the right-hand side, except the first one,
must be of order $n^{-2\alpha}$. This in turn requires
taking $\alpha\leq 1/3$ and $N_n$ of the order $n^{6\alpha}$.
In particular this applies when $\alpha=1/3$, so that this rate
is obtainable but at the cost of quickly increasing $N_n$ at rate $n^2$.
We also notice that when $\alpha=1/4$, to obtain the rate of convergence
$n^{-1/4}\log n$ it is required to take $N_n$ of order $n^{3/2}$,
which is larger than the linear rate
used in Theorem~\ref{Th:unfixed-observations}.

However, a more fair way to look at convergence rates is to express
them in terms of the number of numerical operations performed.
We assume that the computational cost of computing an approximation
$\psi^N(x)$ is of order $N$; this is for instance the case for the
Monte Carlo schemes discussed in Appendix~\ref{app:coupling}.
With $N_n$ being affine in $n$, the total computational cost 
up to stage $n$ of the simulated annealing scheme is thus of order $n^2$.
Denoting the total number of numerical operations performed by $C$,
we then find that the convergence rate is of order $C^{-1/8}\log C$.
If we rather use the second bound above, which requires $N_n$ of order
$n^{6\alpha}$, we see that the computational cost up to stage $n$
is of order $n^{6\alpha+1}$ and the
convergence rate is of order $C^{-\alpha/(6\alpha+1)}\log C$
for $0 <\alpha\leq 1/3$. The optimal $\alpha$ is $\alpha=1/3$, with
rate $C^{-1/9}\log C$. This is inferior to $C^{-1/8}\log C$, so that
the decomposition of the proof of Theorem~\ref{Th:unfixed-observations}
is superior; it does require Hypothesis~\ref{Hyp:bias} however.

\section{A numerical illustration}\label{sec:numerical}

In this section we consider simulated annealing applied
to the likelihood function of a state-space model
as in Appendix~\ref{app:particlefilter}. Thus assume that we
have an observed sequence $(y_t)_{1\leq t\leq T}$ from
a state-space model $((S_t,Y_t))_{1\leq t\leq T}$, whose
Markov transition kernel $Q$ and conditional output densities
$r(\cdot|s)$ both depend on an unknown parameter (vector) $\theta$
which we wish to estimate using maximum likelihood.

The log-likelihood function that we aim to maximise is
$$
\ell_T(\theta) = \sum_{t=1}^T \log p_\theta(y_t|y_{1:t-1})
= \sum_{t=1}^T \log \int r_\theta(y_t|s) \, \pi_{t|t-1}^{\theta}(ds),
$$
where $p_\theta(y_t|y_{1:t-1})$ is the conditional density of
$Y_t$ given $Y_{1:t-1}$, and $\pi_{t|t-1}^\theta$ is
the predictive distribution $\p_\theta(S_t\in\cdot\mid y_{1:t-1})$.
As $\pi_{t|t-1}^\theta$ can in general not be computed we
need to approximate the log-likelihood function, and one
way to do that is through
$$
\ell_T^N(\theta) 
= \sum_{t=1}^T \log \int r_\theta(y_t|s) \, \pi_{t|t-1}^{\theta,N}(ds),
$$
where we take $\pi_{t|t-1}^{\theta,N}(ds)$ as the particle
filter approximation of Appendix~\ref{app:particlefilter}.

The log-likelihood function is essentially a sum of functions of
the form studied in Appendix~\ref{app:particlefilter}, except for
the logarithmic transformation. Assuming however,
as in Appendix~\ref{app:particlefilter}, that $r_\theta$ is
uniformly bounded from below by some $\underline{r}>0$, we find
that each of the integrals above are bounded from below by
$\underline{r}$. Moreover, using the inequality
$|\!\log x-\log y|\leq |x-y|/(x\wedge y)$, valid for all $x,y>0$,
we find that
$$
|\ell_T^N(\theta) - \ell_T^{N'}(\theta)|
\leq \frac{1}{\underline{r}} \sum_{t=1}^T
\left| \int r_\theta(y_t|s) \, \pi_{t|t-1}^{\theta,N}(ds)
  - \int r_\theta(y_t|s) \, \pi_{t|t-1}^{\theta,N'}(ds) \right|.
$$
This bound involves sums of functions of the form studied in
Appendix~\ref{app:particlefilter}
(take $h(s)=r_\theta(y_t|s)$), and we can proceed as there to
show that Hypothesis~\ref{Hyp:bias} holds.
A similar argument where we replace $\ell_T(\theta)$ by the
exact likelihood and appeal to
Theorem~7.4.4 of \citet{delmoral-04} shows that 
Hypothesis~\ref{Hyp:approx} holds.

\subsection{Simulation study}

We considered the benchmark model
\citep[Eqs.~8.3.4--8.3.5]{doucet-freitas-gordon-01}
\begin{eqnarray}
S_t & =  & a S_{t-1} + b\frac{S_{t-1}}{1+S_{t-1}^2}
     +\gamma\cos(1.2t)+\sigma_v V_t,  \label{eq:sysdyn} \\
Y_t & = & \frac{S_t^2}{20}+\sigma_w W_t, \label{eq:obseq}
\end{eqnarray}
where $(S_t)$ is the unobserved Markov chain taking
values in $\R$, $(Y_t)$ is the observable process and
$(V_t)$ and $(W_t)$ are mutually independent sequences of
i.i.d.\ standard Gaussian random variables. 
We wish to estimate the five model parameters
$\theta=(a,b,\gamma,\sigma_v,\sigma_w)$ given a sequence
$(y_t)_{1\leq t\leq T}$ of observations, and we did so using
the approximate maximum likelihood (ML) approach outlined above
with the bootstrap particle filter, i.e.\ particle mutations following
the system dynamics (\ref{eq:sysdyn}).
We remark that the state space of the model above is not compact, 
so that the conditional densities $r_\theta(y|s)$ are not bounded from
below in $s$. The model does thus not fulfil the technical
conditions made above, but the results below are still
an illustrative example of how the simulated annealing
scheme performs in a particular case.

\begin{figure}[htb]
\begin{center}
\resizebox{80mm}{!}{\centerline{\includegraphics{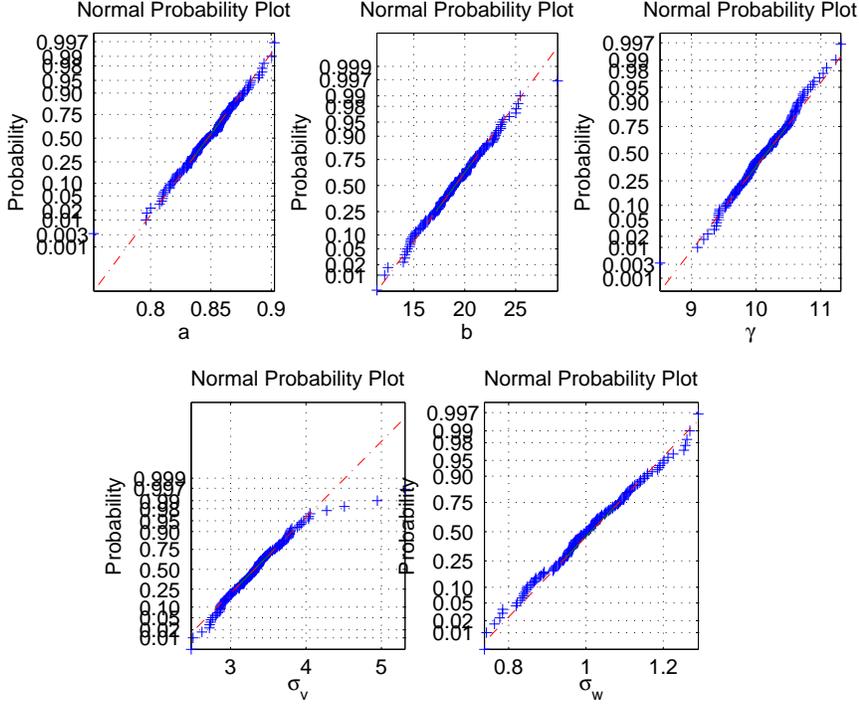}}}
\end{center}
\caption{Normal probability plots of approximate ML estimates of
parameters $(a,b,\gamma,\sigma_v,\sigma_w)$ in the model
(\ref{eq:sysdyn})--(\ref{eq:obseq}), obtained from 150 replications
of 5,000 iterations of the simulated annealing scheme applied
to the particle filter approximation of the log-likelihood.}
\label{Fig:Normplotpar}
\end{figure}

We simulated a single trajectory $(y_t)_{1\leq t\leq T}$ of length
$T=500$ with parameters 
$\theta^0=(a^0,b^0,\gamma^0,\sigma_v^0,\sigma_w^0)
=(0.9, 18, 10, \sqrt{10}, 1)$.
In the simulated annealing scheme we let the inverse temperature be
$\beta_n = 10 n^{1/4}$, corresponding to $\tau=1/10$ in 
Corollary~\ref{Cor:optimal}, and let number of particles at step $n$ be
$N_n = n\vee 20$, a function which is affine for $n\geq 20$
(Theorem~\ref{Th:unfixed-observations}). The algorithm was run
for 5,000 iterations in each of 150 independent replications.
The parameter space $\Theta$ was taken as the five-dimensional
hyper-rectangle 
$[0.45, 1.8]\times[9, 36]\times[5, 20]\times[0.316, 36]\times [0.5, 2]$.
For $K$ we used a Gaussian random walk proposal (on the log-scale
for the standard deviations), where we constrained the random walk to
$\Theta$; any coordinate of the parameter proposed outside
$\Theta$ was pulled back to the boundary.
The incremental covariance of the kernel at step $n$ was a diagonal matrix
whose $i$-th diagonal element was the squared $i$-th side length of
$\Theta$ divided by $\log(n+1)^2$. In each replication the initial
point $\theta_0$ was drawn uniformly on $\Theta$.

After 5,000 iterations of the simulated annealing algorithm,
the sample means and standard errors of the parameter estimates
$\bt_{5000}$ (over the 150 replications) were
$(0.85, 19.1, 10.1, 3.4, 1.01)$ and $(0.024, 3.0, 0.46, 0.41, 0.11)$
respectively. These sample means are in good
agreement with the true $\theta^0$. Ideally we would like to compare
to the ML estimates, which are however unavailable.
Figure~\ref{Fig:Normplotpar} shows that the estimates follow normal
distributions with good accuracy, with the exception of $\sigma_v$.
This of course is an empirical observation for which we have
no theoretical support, as we have not discussed
convergence in law of the differences $\theta_n-\theta_{\max}$
and $\bt_n-\theta_{\max}$, suitably scaled, where $\theta_{\max}$
is the point where $\psi$ is maximal.



\bibliographystyle{rss_ours}
\bibliography{bibliography}

\appendix

\section{Rate of convergence of classical simulated annealing}
\label{app:classicalrate}

In this section we prove the bound (\ref{eq:stdrateconv}) and also,
by studying a specific example, that this bound cannot be improved
generally. We assume that Hypotheses~\ref{Hyp:mixing}--\ref{Hyp:osc}
and Hypothesis~\ref{Hyp:gradpsi}(ii) hold. Since we now consider
classical simulated annealing we have $f(t)=\exp(-t)$, and we
take $\beta_n = \beta_0 \log(n+e)$ with $1/\beta_0>\osc(\psi)$
\citep[cf.\ ][Theorem~2.3.5]{bartoli-delmoral-01}.
As in Section~\ref{sec:rate} we let $\eta_n$ be the law of 
$\theta_n$ and denote by $\mu_\beta$ the invariant distribution of $K_\beta$.

Now write
$$
\p(\psi(\theta_n)\leq \psimax -\epsilon)
= \eta_n(\uec) = ( \eta_n(\uec) - \mu_{\beta_n}(\uec) )
   + \mu_{\beta_n}(\uec).
$$
We will show that the first term of this decomposition (the difference)
tends to zero at algebraic rate, while the second term vanishes only
logarithmically fast. Thus the left-hand side tends to zero
at logarithmic rate too. In a specific example we will also show that the
logarithmic rate for the second term, which in general is a bound,
is in fact the exact rate; thus the logarithmic rate for the left-hand
side cannot be improved generally.
Here  emerges an essential difference between classical simulated
annealing and the new scheme analysed in Section~\ref{sec:rate}.
In both cases the total variation distance between
the law $\eta_n$ of $\theta_n$ and the invariant law
$\mu_{\beta_n}$ vanishes at algebraic rate; $n^{\alpha-1}$ 
for classical simulated annealing (see below) and $n^{2\alpha-1}$
for the new scheme (Theorem~\ref{Th:rate}).
The rate at which $\mu_{\beta_n}$ concentrates around the
maximum of $\psi$ is much different however; this rate is
algebraic too for the new scheme (Lemma~\ref{Lem:stationary}),
but only logarithmic (or algebraic with rate tending to zero)
for the classical scheme.

We now proceed to the details.
Put once again $\Delta_n=\|  \eta_n - \mu_{\beta_n} \|_{\TV}$.
We then have the recursion
$$
\Delta_{n+1}
\leq (1-\epsilon_K e^{-\beta_n \osc(\psi)}) \Delta_n
+ (\beta_{n+1}-\beta_n) \osc(\psi);
$$
see \citet[Remark~3.3.13]{bartoli-delmoral-01} and cf.\ the proof
of Theorem~\ref{Th:rate}. With  the present choice of $(\beta_n)$
we find $\beta_{n+1}-\beta_n\leq \beta_0/(n+1)$ and
$\exp(-\beta_n\osc(\psi))=(n+e)^{-\alpha}$ with 
$\alpha=\beta_0\osc(\psi)<1$. Iterating the above recursion yields
\begin{eqnarray*}
\Delta_{n+1} & \leq &
  \sum_{q=1}^{n} \prod_{k=q+1}^n
  \left(1-\frac{\epsilon_K}{(k+e)^\alpha} \right)
  \times \frac{\beta_0\osc(\psi)}{q+1} \\
& + & \prod_{k=1}^n \left(1-\frac{\epsilon_K}{(k+e)^\alpha} \right)
  \times \|\eta_1 - \mu_{\beta_1}\|_{\TV},
\end{eqnarray*}
where an empty product (when $q=n$) is interpreted as unity.
Bound the product as
\begin{eqnarray*}
\log \prod_{k=q+1}^n \left(1-\frac{\epsilon_K}{(k+e)^\alpha} \right)
& \leq & -\sum_{k=q+1}^{n} \frac{\epsilon_K}{(k+e)^\alpha} \\
& \leq & -\epsilon_K \int_{q+1}^{n+1} \frac{dx}{(x+e)^\alpha} \\
& = & -C_1 ((n+e+1)^{1-\alpha} - (q+e+1)^{1-\alpha}),
\end{eqnarray*}
where $C_1=\epsilon_K/(1-\alpha)$. Thus, using $\beta_0\osc(\psi)<1$ again
as well,
\begin{eqnarray*}
\Delta_{n+1} & \leq & e^{-C_1(n+e+1)^{1-\alpha}}
  \sum_{q=1}^n e^{C_1(q+e+1)^{1-\alpha}} \frac{1}{q+1} \\
  & + & e^{-C_1((n+e+1)^{1-\alpha} - (e+1)^{1-\alpha})}
  \|\eta_1 - \mu_{\beta_1}\|_{\TV} \\
& \leq & e^{-C_1(n+e+1)^{1-\alpha}}
   \int_1^{n+1} e^{C_1(x+e+1)^{1-\alpha}} \frac{1}{x} \, dx
   + C e^{-C_1(n+e+1)^{1-\alpha}}.
\end{eqnarray*}
By manipulating the integral on the right-hand side, $I_{n+1}$ say,
we can just as in the proof of Theorem~\ref{Th:rate} prove that
$$
I_{n+1} \leq C e^{C_1(n+e+2)^{1-\alpha}} \frac{1}{(n+1)^{1-\alpha}}.
$$
Hence we obtain
\begin{eqnarray*}
\Delta_{n+1} & \leq & C e^{-C_1((n+e+1)^{1-\alpha} - (n+e+2)^{1-\alpha})}
    \frac{1}{(n+1)^{1-\alpha}} + C e^{-C_1(n+e+1)^{1-\alpha}}\\
& \leq & \frac{C}{(n+1)^{1-\alpha}}
\end{eqnarray*}
and thus $\Delta_n\leq C/n^{1-\alpha}$.

So far the difference between $\eta_n$ and $\mu_{\beta_n}$. We now
turn to how concentrated $\mu_{\beta_n}$ is around the maximum
of $\psi$. To start with we may employ Lemma~\ref{Lem:stationary},
with $f$ and $\beta_n$ as above, to obtain
$$
\mu_{\beta_n}(\uec) \leq \frac{C_\epsilon}{\beta_0\log(n+e)}
  + (n+e)^{-\beta_0(\epsilon''-\epsilon)};
$$
a logarithmic rate in other words. We can also use the property
mentioned in \cite[p.~64]{bartoli-delmoral-01}, that $\mu_\beta$ equals
$\exp(\beta\psi(x))\,\gamma(dx)$ up to a normalising constant
with $\gamma$ the invariant distribution of $K$, to obtain
\begin{eqnarray}
\mu_\beta(\uec) & = & \frac{\ds\int_{\uec}
    e^{\beta\psi(y)} \,\gamma(dy)}
    {\ds \int e^{\beta\psi(y)}\,\gamma(dy)}
    \nonumber \\[3mm]
&  \leq & \frac{\ds\int_{\uec}
    e^{\beta(\psi_{\max}-\epsilon)}\,\gamma(dy)}
    {\ds \int_{U^{\epsilon/2}}
        e^{\beta(\psi_{\max}-\epsilon/2)}
        \,\gamma(dy)}
 = \frac{e^{-\beta\epsilon/2}}{\gamma(U^{\epsilon/2})}.
 \label{eq:classicalmubound}
\end{eqnarray}
Inserting $\beta_n$ for $\beta$, it follows that
$$
\mu_{\beta_n}(\uec) \leq \frac{1}{\gamma(U^{\epsilon/2})}
  (n+e)^{-(\beta_0/2)\epsilon}.
$$
This is the bound (\ref{eq:stdrateconv}).

We now prove that this bound cannot be improved in general. Consider
the example $\Theta=[-1/2,\,1/2]$, $\psi(x)=-|x|$,
$K(x,dy)=dy$. Thus $K$ is an independence kernel that proposes
uniformly on $\Theta$. It is immediate that the invariant measure
$\gamma$ of $K$ is Lebesgue measure on $\Theta$, and that $\gamma$
is $K$-reversible. Now $\mu_\beta(A)$ is proportional to
$\int_A \exp(-\beta|y|)\,dy$, so that
$$
\mu_\beta(\uec) = \frac
   {\ds \int_{\epsilon<|y|\leq 1/2} e^{-\beta |y|}\, dy}
   {\ds \int_\Theta e^{-\beta |y|}\, dy}
= \frac{e^{-\beta\epsilon}-e^{-\beta/2}}{1-e^{-\beta/2}}
\sim e^{-\beta \epsilon} \quad\mbox{as $\beta\to\infty$}.
$$
We can indeed, by an obvious modification of the argument above,
adjust (\ref{eq:classicalmubound}) into the bound
$1/\gamma(U^{\epsilon\delta}) \times e^{-\beta(1-\delta)\epsilon}$,
where $0<\delta<1$ is arbitrary. The rate of this bound thus
can thus be made arbitrarily close to the exact rate of this example.

\section{Coupling function approximations}\label{app:coupling}

The purpose of this appendix is to illustrate how one may construct
function approximations $\psi^N$ that satisfy Hypothesis~\ref{Hyp:bias},
and how the relatively `high level' condition of this hypothesis can
be guaranteed by more `low level' assumptions.

Thus assume that we are given a probability measure $\mu$ on $\Theta$,
$\beta>0$, and two approximation indices $N$ and $N'$. We wish to bound
$\|\mu K_\beta^N-\mu K_\beta^{N'}\|_{\TV}
=\sup_A |\mu K_\beta^N(A)-\mu K_\beta^{N'}(A)|$, where the supremum is over
$A\in\mathcal{B}(\Theta)$. We will accomplish this by
constructing two coupled samples from $\mu K_\beta^N$ and
$\mu K_\beta^{N'}$ respectively as follows.
\begin{enumerate}
  \item[(i)] Sample a point $x$ from $\mu$ and then a point $z$
    from $K(x,\cdot)$.
  \item[(ii)] Compute the function approximations $\psi^N(x)$,
    $\psi^N(z)$, $\psi^{N'}(x)$ and $\psi^{N'}(z)$. For the time being
    we do not specify exactly how this is done.
  \item[(iii)] Sample a random number $U$ from the uniform distribution
    on $(0,1)$ and accept the proposal $z$ if
    $U\leq f(\beta(\psi^N(x)-\psi^N(z))_+)$ or
    $U\leq f(\beta(\psi^{N'}(x)-\psi^{N'}(z))_+)$ respectively,
    for the two indices $N$ and $N'$.
\end{enumerate}
The samples $\mu K_\beta^N$ and $\mu K_\beta^{N'}$ so
constructed will be different only if the two decisions is step~(iii)
are different, so the probability of the former event is bounded by the
probability of the latter one.
To compute the probability that the decisions of step~(iii) differ,
we notice this event occurs if $U$ falls in between
the two function values used there, which, since $U$ is uniform,
happens with (conditional) probability
$$
|f(\beta(\psi^N(x)-\psi^N(z))_+) - f(\beta(\psi^{N'}(x)-\psi^{N'}(z))_+)|.
$$
Hence the probability of different decisions in step~(iii) is bounded by
$$
\sup_{x,z\in\Theta}
\E|f(\beta(\psi^N(x)-\psi^N(z))_+) - f(\beta(\psi^{N'}(x)-\psi^{N'}(z))_+)|,
$$
where the expectation is w.r.t.\ the function approximations
$\psi^N$ and $\psi^{N'}$.
The difference of the function values can be bounded as
$$
\beta |(\psi^N(x)-\psi^N(z))_+ - (\psi^{N'}(x)-\psi^{N'}(z))_+|
 \times f'(\zeta),
$$
where $\zeta$ is point between the two function arguments. By the
assumptions on $f$ its derivative is necessarily bounded,
and it is straightforward to check that for any real $a$ and $b$,
$|a_+-b_+|\leq |a-b|$. Therefore the probability of different
decisions in step~(iii) is bounded by
\begin{eqnarray*}
\lefteqn{ \beta \| f'\|_\infty
 \sup_{x,z\in\Theta}
 \E |(\psi^N(x)-\psi^N(z)) - (\psi^{N'}(x)-\psi^{N'}(z))| }
 \hspace*{30mm} \\
& \leq  & 2 \beta \| f'\|_\infty \sup_{x\in\Theta}
 \E |\psi^N(x)-\psi^{N'}(x)|.
\end{eqnarray*}
Thus, at this point we see that if the function approximations satisfy
\begin{equation}\label{eq:couplingcondition}
\sup_{x\in\Theta}\E |\psi^N(x)-\psi^{N'}(x)|
\leq C \frac{N'-N}{N}
\end{equation}
for some constant $C$, Hypothesis~\ref{Hyp:bias} will hold.

Verifying (\ref{eq:couplingcondition}) is, of
course, a problem very much related to the specific construction
of these approximations. In the following two subsections we will
deal with two specific settings: i.i.d.\ sample means and
particle filters.

\subsection{Simple Monte Carlo sample means}

Here we consider the possibly simplest of all approximation
schemes: a sample mean of i.i.d.\ summands. Thus we assume
that for a random variable $\xi$ with some known distribution
and some known function $h$,
$\psi(x)=\E h(\xi;x)$ where the expectation is w.r.t.\ $\xi$,
and that its approximation is
$$
\psi^N(x) = \frac{1}{N} \sum_{i=1}^N h(\xi_i;x)
$$
where the $\xi_i$ are i.i.d.\ variables distributed as $\xi$.
We note in passing that for this scheme the Marcinkiewicz-Zygmund
inequality \citep[p.~498]{shiryaev-1995} with $p=1$
implies that Hypothesis~\ref{Hyp:approx} holds.
Moreover, for $N'>N$,
$$
\psi^N(x) - \psi^{N'}(x)
= \left(\frac{1}{N}-\frac{1}{N'}\right)\sum_{i=1}^N h(\xi_i;x)
- \frac{1}{N'}\sum_{i=N+1}^{N'} h(\xi_i;x)
$$
and
\begin{eqnarray*}
  \E|\psi^N(x) - \psi^{N'}(x)|
  & \leq & \left(\frac{1}{N}-\frac{1}{N'}\right) N \E|h(\xi;x)|
  + \frac{1}{N'}(N'-N) \E|h(\xi;x)| \\
& = & 2\,\E|h(\xi;x)| \, \frac{N'-N}{N'} .
\end{eqnarray*}
It is now immediate that if $\E|h(\xi;x)|$ is bounded in
$x\in\Theta$, (\ref{eq:couplingcondition}) holds.

\subsection{Particle filter estimates}\label{app:particlefilter}

Consider a state-space model $((S_t,Y_t))_{t\geq 1}$, where
$(S_t)$ is an unobserved Markov chain on some general state space and
$(Y_t)$ is an observed sequence of random variables. The
association between $(S_t)$ and $(Y_t)$ is local in the sense that
(i) given $(S_t)$, the $Y$-variables are conditionally independent,
and (ii) given $(S_t)$ and for any time index $u$, the
conditional distribution of $Y_u$ depends on $S_u$ only.

We will denote the transition kernel of the Markov chain $(S_t)$
by $Q$, and the conditional density of $Y_t$ given $S_t=s$
by $r(\cdot|s)$. Both of these quantities are assumed to depend
on some model parameters $\theta$, which we indicate by writing
$Q_\theta$ and $r_\theta$ respectively.

The function $\psi$ we wish to approximate
is $\psi(\theta) = \E_\theta[h(S_t) \mid y_{1:t-1}]$, that is,
the expectation of some function $h$ w.r.t.\ the
so-called \textit{predictive distribution}
$\pi_{t|t-1}^\theta(\cdot)=\p_\theta(S_t\in\cdot\mid y_{1:t-1})$,
where $t\geq 1$ is some time index, the notation $y_{1:t-1}$
is short for $y_1,y_2,\ldots,y_{t-1}$, and subindex `$t|t-1$'
indicates that the distribution concerns the state at time $t$
conditional on observed data up to time $t-1$.

The predictive distributions can,
together with the so-called \textit{filter distributions}
$\pi_{t|t}^\theta(\cdot)=\p_\theta(S_t\in\cdot\mid y_{1:t})$,
be computed recursively in time---at least in principle.
The recursive formulae read
\begin{equation}\label{eq:pred2filt}
  \pi_{t|t}^\theta(ds)
     = \frac{r_\theta(y|s) \, \pi_{t|t-1}^\theta(ds)}
     {\ds \int r_\theta(y|s') \, \pi_{t|t-1}^\theta(ds')}
\end{equation}
and
\begin{equation}\label{eq:filt2pred}
\pi_{t+1|t}^\theta(\cdot)
= \int Q_\theta(s,\cdot) \, \pi_{t|t}^\theta(\cdot).
\end{equation}
The first of these formulae is just Bayes' rule, and the second
one means to propagate the filter through the state dynamics
$Q_\theta$.

In practice the above relations do no admit exact numerical solution
except in two cases: when the state space of $(S_t)$ is finite
(so-called hidden Markov models; the integrals then turn into finite
sums) and when the state-space model is linear with additive
Gaussian noise (the solution then being provided by the Kalman filter).
There are many ways to approximate these two recursions, and
here we shall examine an approach referred to as \textit{particle filters}.
This section contains a full introduction neither to state-space models
nor to particle filters, and we refer to \cite{doucet-freitas-gordon-01}
for a more complete coverage of both.

The basic idea of a particle filter is to approximate the filter
and predictive distributions with the empirical distributions of
a set of \textit{particles}, whose positions are dynamically
updated in time. There is not just one particle filter algorithm---the
term rather refers to a framework for algorithms---and the particular
algorithm we look at here is usually denoted the
\textit{bootstrap particle filter}. We now describe how this algorithm
works; the parameter $\theta$ and population size $N$ are fixed
throughout.

Assume that at some time index $t$ we have available a collection\\
$(\xi_{t|t-1,i}^{\theta,N})_{1\leq i\leq N}$ of particles whose
empirical distribution approximates $\pi_{t|t-1}^\theta$.
The transformation (\ref{eq:pred2filt}) is approximated as follows.

\begin{itemize}
  \item[(a)] \textit{Weighting}. Compute unnormalised weights
    $\tilde{w}_{t,i}^{\theta,N}=r_\theta(y_t|\xi_{t|t-1,i}^{\theta,N})$
    and then normalised weights
    $w_{t,i}^{\theta,N}=\tilde{w}_{t,i}^{\theta,N}/
    \sum_j \tilde{w}_{t,j}^{\theta,N}$.
  \item[(b)] \textit{Resampling}. Create a sample
   $(\xi_{t|t,i}^{\theta,N})_{1\leq i\leq N}$ by sampling
   $N$ times independently from $(\xi_{t|t-1,i}^{\theta,N})_{1\leq i\leq N}$
   with weights $(w_{t,i}^{\theta,N})_{1\leq i\leq N}$.
\end{itemize}
The empirical distribution of the sample
$(\xi_{t|t,i}^{\theta,N})_{1\leq i\leq N}$ obtained in the resampling
step approximates $\pi_{t|t}^\theta$.

The transformation (\ref{eq:filt2pred}) is approximated as
follows.
\begin{itemize}
  \item[(c)] \textit{Mutation}. Create a sample
   $(\xi_{t+1|t,i}^{\theta,N})_{1\leq i\leq N}$ by independently sampling
   $\xi_{t+1|t,i}^{\theta,N}$ from
   $Q_\theta(\xi_{t|t,i}^{\theta,N},\cdot)$.
\end{itemize}

The procedure is initialised at time $t=0$ by letting
$(\xi_{1|0,i}^{\theta,N})_{1\leq i\leq N}$ be an i.i.d.\ sample of
size $N$ from the initial distribution $P_\theta(S_1\in\cdot)$
of the state process. This distribution may depend on $\theta$ but
is otherwise assumed known.

The book by \citet{delmoral-04} is a thorough treatise of theoretical
properties of particle filters, and in particular its Theorem~7.4.4
shows that Hypothesis~\ref{Hyp:approx} holds, provided
that for each $y_t$, $r_\theta(y_t|s)$ is bounded in $\theta$ and $s$.
We are here particularly interested in the particle approximations
of the predictive distributions, and the update of these can be summarised
as follows: compute the normalised weights $w_{t,i}^{\theta,N}$
and then sample for $1\leq i\leq N$, independently, first an index
$j$ with probability $w_{t,j}^{\theta,N}$ and then
$\xi_{t+1|t,i}^{\theta,N}\sim Q_\theta(\xi_{t|t-1,j}^{\theta,N})$.

We will now run, simultaneously, two particles filters of sizes
$N'>N$ respectively. All other
properties of the filers---data, parameters, dynamics---agree.
The joint dynamics of the filters will be coupled in a way
such that many particles of the two filters, at any given time index,
coincide. Indeed, for each time index $t$ we define a partition
$J_t\cup J_t^c$ of $\{1,2,\ldots,N'\}$ such that
$\xi_{t|t-1,i}^{\theta,N} = \xi_{t|t-1,i}^{\theta,N'}$ for
$i\in J_t$. The details of the coupling are as follows.

\begin{itemize}
  \item[(i)] \textit{Initialisation}. Sample
    $(\xi_{1|0,i}^{\theta,N'})_{1\leq i\leq N'}$ independently
    from $\p_\theta(S_1\in\cdot)$, let
    $\xi_{1|0,i}^{\theta,N} = \xi_{1|0,i}^{\theta,N'}$ for
    $1\leq i\leq N$ and let $J_1=\{1,2,\ldots,N\}$,
    $J_1^c=\{N+1,N+2,\ldots,N'\}$.

  \item[(ii)] \textit{Recursion from $t$ to $t+1$}.
    We have $\xi_{t|t-1,i}^{\theta,N} = \xi_{t|t-1,i}^{\theta,N'}$ for
    $i\in J_t$ and compute the weights
    $(w_{t,i}^{\theta,N})_{1\leq i\leq N}$ and
    $(w_{t,i}^{\theta,N'})_{1\leq i\leq N'}$.

    When sampling the new particles, we couple the two filters
    in a way such that independently for each $1\leq i\leq N$, 
    one of the events below take place
    (index  $j$ has the same meaning as above):
    \begin{itemize}
      \item for $j\in J_t$,
        \begin{enumerate}
        \item[--] $\xi_{t+1|t,i}^{\theta,N} = \xi_{t+1|t,i}^{\theta,N'}
          \sim Q_\theta(\xi_{t|t-1,j}^{\theta,N},\cdot)$
          with probability  $w_{t,j}^{\theta,N} \wedge w_{t,j}^{\theta,N'}$;
        \item[--] $\xi_{t+1|t,i}^{\theta,N} 
          \sim Q_\theta(\xi_{t|t-1,j}^{\theta,N},\cdot)$
          with probability
          $w_{t,j}^{\theta,N}-w_{t,j}^{\theta,N} \wedge w_{t,j}^{\theta,N'}$;
        \item[--] $\xi_{t+1|t,i}^{\theta,N'} 
          \sim Q_\theta(\xi_{t|t-1,j}^{\theta,N'},\cdot)$
          with probability
          $w_{t,j}^{\theta,N'}-w_{t,j}^{\theta,N} \wedge w_{t,j}^{\theta,N'}$;
        \end{enumerate}
      \item for $j\in J_t^c\cap\{1,2,\ldots,N\}$,
        \begin{enumerate}
        \item[--] $\xi_{t+1|t,i}^{\theta,N} 
          \sim Q_\theta(\xi_{t|t-1,j}^{\theta,N},\cdot)$
          with probability  $w_{t,j}^{\theta,N}$;
        \item[--] $\xi_{t+1|t,i}^{\theta,N'} 
          \sim Q_\theta(\xi_{t|t-1,j}^{\theta,N'},\cdot)$
          with probability  $w_{t,j}^{\theta,N'}$.
        \end{enumerate}
     \end{itemize}
     Finally, for $N<i\leq N'$,
     $\xi_{t+1|t,i}^{\theta,N'} \sim Q_\theta(\xi_{t|t-1,j}^{\theta,N'},.)$
     with probability $w_{t,j}^{\theta,N'}$.

     We let $J_{t+1}$ be the set of indices $1\leq i\leq N$
     such that the first of the above events happened.
\end{itemize}

From this construction it is immediate that the distributions of
the two filters are the same as if they had been run separately
and independently in the usual manner. Let
$\pi_{t|t-1}^{\theta,N}$ be the particle filter
approximation to the predictive distribution at time index $t$;
$$
\pi_{t|t-1}^{\theta,N}(A) =
  \frac{1}{N} \sum_{i=1}^N I_A(\xi_{t|t-1,i}^{\theta,N})
$$
for all $A\in\mathcal{B}(\Theta)$, where $I_A$ is the
indicator function of $A$.

\begin{proposition}\label{prop:pfcoupling}
Assume that observations $y_{1:T}$ are given and that there
is a number $\underline{r}>0$ such that
$\underline{r}\leq r_\theta(y_t|s)\leq 1/\underline{r}$
for all $1\leq t\leq T$, all $s$ in the state space and all $\theta\in\Theta$.
Then there are constants $C_t$ for $1\leq t\leq T$ such that
for any integers $N'>N>0$,
$$
\E\| \pi_{t|t-1}^{\theta,N} - \pi_{t|t-1}^{\theta,N'}\|_{\TV}
\leq C_t \left( \frac{N'-N}{N} \right).
$$
\end{proposition}
The constants $C_t$ depend on $\underline{r}$, but otherwise the
bound is uniform in $\theta$. Therefore this result implies
(\ref{eq:couplingcondition}) for
$\psi^N(\theta) = \int h(s) \, \pi_{t|t-1}^{\theta,N}(ds)$
whenever $h$ is bounded on the state space of $(S_t)$.

The requirement of a lower bound $\underline{r}>0$ on $r_\theta$,
uniform in $\theta$ and $s$, will typically be satisfied only
if both $\Theta$ and the state space of $(S_t)$ are compact,
or at least bounded. Boundedness of $\Theta$ is as good as implied
by Hypothesis~\ref{Hyp:mixing}, whereas boundedness of the
state space is a more serious limitation. Having said that we notice that
this condition is recurring in the literature on particle filters,
in particular when treating forgetting properties.

\begin{proof}[Proof of Proposition \ref{prop:pfcoupling}]
For any $A\in\mathcal{B}(\Theta)$,
\begin{eqnarray}
  |\pi_{t|t-1}^{\theta,N}(A) - \pi_{t|t-1}^{\theta,N'}(A)|
  & \leq &
  \left| \frac{1}{N} \sum_{i=1}^N I_A(\xi_{t|t-1}^{\theta,N})
     - \frac{1}{N'}
     \sum_{i=1}^{N'} I_A(\xi_{t|t-1}^{\theta,N'}) \right| \nonumber \\
  &= & \Bigg| \sum_{i\in J_t}^N I_A(\xi_{t|t-1}^{\theta,N})
  \left( \frac{1}{N}-\frac{1}{N'} \right) \nonumber \\
  && + \frac{1}{N} \sum_{i\in J_t^c,\,i\leq N} I_A(\xi_{t|t-1}^{\theta,N})
   + \frac{1}{N'} \sum_{i\in J_t^c} I_A(\xi_{t|t-1}^{\theta,N'})
   \Bigg| \nonumber \\
   &  \leq & \# J_t \left( \frac{1}{N}-\frac{1}{N'} \right)
     + \# J_t^c \left( \frac{1}{N} + \frac{1}{N'} \right) \nonumber \\
   & \leq & 1 - \frac{N}{N'}
     + \# J_t^c \left( \frac{1}{N} + \frac{1}{N'} \right),
     \label{eq:pftvbound1}
\end{eqnarray}
where $\#$ denotes cardinality of a set and,
in the last step, $\# J_t$ was bounded by $N$. We now seek to bound
$\E(\#J_t^c)$.

Put $p_t = 1-\#J_t/N'$ and define the $\sigma$-field
$\mathcal{F}_t = \sigma(\xi_{t|t-1,i}^{\theta,N},\, 1\leq i\leq N)
\vee \sigma(\xi_{t|t-1,i}^{\theta,N'},\, 1\leq i\leq N')$.
Then conditionally on $\mathcal{F}_t$, $\# J_{t+1}$ is a binomial random variable
with parameters $N$ and
$\sum_{i\in J_t} (w_{t,i}^{\theta,N} \wedge w_{t,i}^{\theta,N'})$.
Using the definition of $J_t$ and abbreviating
$r_\theta(y_t|s)$ as $r_{t,\theta}(s)$, we find that
\begin{eqnarray*}
\sum_{i\in J_t} w_{t,i}^{\theta,N} \wedge w_{t,i}^{\theta,N'}
& = & \frac{\sum_{i\in J_t} r_{\theta,t}(\xi_{t|t-1,i}^{\theta,N})}
  {\sum_{1\leq i\leq N} r_{\theta,t}(\xi_{t|t-1,i}^{\theta,N})}
  \left( 1 \wedge
    \frac{\sum_{1\leq i\leq N} r_{\theta,t}(\xi_{t|t-1,i}^{\theta,N})}
  {\sum_{1\leq i\leq N'} r_{\theta,t}(\xi_{t|t-1,i}^{\theta,N'})}
  \right) \\
& = &\left(  \frac{1}
   {1 + \frac{\sum_{i\in J_t^c,\,i\leq N}
              r_{\theta,t}(\xi_{t|t-1,i}^{\theta,N})}
  {\sum_{i\in J_t} r_{\theta,t}(\xi_{t|t-1,i}^{\theta,N})} }\right) \\
&&~~\times
  \left( 1 \wedge
    \frac{1 +
      \frac{\sum_{i\in J_t^c,\,i\leq N} r_{\theta,t}(\xi_{t|t-1,i}^{\theta,N})}
  {\sum_{i\in J_t} r_{\theta,t}(\xi_{t|t-1,i}^{\theta,N})} }
  {1 + \frac{\sum_{i\in J_t^c} r_{\theta,t}(\xi_{t|t-1,i}^{\theta,N'})}
  {\sum_{i\in J_t} r_{\theta,t}(\xi_{t|t-1,i}^{\theta,N})} } \right) \\
& \geq & \frac{1}{1 + \underline{r}^{-2} \frac{p_t}{1-p_t}}
   \times \frac{1}{1 + \underline{r}^{-2} \frac{p_t}{1-p_t}} \\
& = & \frac{1}{\left(1 + \underline{r}^{-2} \frac{p_t}{1-p_t}\right)^2}
=: u(p_t).
\end{eqnarray*}
We note that as $u$ is convex and decreasing with $u(0)=1$,
there exists a constant $C_u>0$ such that $u(p)\geq 1-C_u p$.

The above-mentioned conditional binomial distribution of $\# J_t$
implies, together with the above inequality,
that $\E(\# J_{t+1}\mid\mathcal{F}_t)\geq N u(p_t)$, and therefore
\begin{eqnarray*}
  \E(p_{t+1} \mid \mathcal{F}_t) & \leq &
  1 - \frac{N u(p_t)}{N'} \\
& = & 1 - u(p_t) + \frac{N'-N}{N'} u(p_t) \\
& \leq & 1 - u(p_t) + \frac{N'-N}{N'} \\
& =: & v(p).
\end{eqnarray*}
Applying this inequality recursively, it follows that
$\E(p_t)\leq v^{\circ t}(p_0)$, where superindex `$\circ t$'
means $t$-fold function composition.

We notice that $p_0=(N'-N)/N'$ and $v(p)\leq C_u p + p_0 \leq C_v(p+p_0)$
for some constant $C_v$, and by induction we find that there is
a constant $C_{v,t}$ such that
$v^{\circ t}(p)\leq C_{v,t}(p+p_0)$.
Thus $\E (\# J_t^c) = N' \E(p_t)
\leq 2 N'  C_{v,t} p_0 = 2 C_{v,t} (N'-N)$.
The proof is finished by
inserting this bound into the right-hand side of
(\ref{eq:pftvbound1}), then taking the supremum over
$A\in\mathcal{P}(\Theta)$ on the left-hand side and finally
the expectation.

\end{proof}
  
\end{document}